\documentclass[11pt]{article}
\usepackage{amsfonts,amssymb,array,amsmath} 
\newtheorem{Lemma}{Lemma}[section]\newcommand{\bel}{\begin{Lemma}}\newcommand{\eel}{\end{Lemma}}
\newtheorem{Remark}[Lemma]{Remark}\newcommand{\beR}{\begin{Remark}\rm}\newcommand{\eeR}{\end{Remark}}
\newtheorem{Definition}[Lemma]{Definition}\newcommand{\bed}{\begin{Definition}\rm}\newcommand{\eed}{\end{Definition}}
\newtheorem{Claim}[Lemma]{Proposition}\newcommand{\bec}{\begin{Claim}}\newcommand{\eec}{\end{Claim}}
\newtheorem{Example}[Lemma]{Example}\newcommand{\bex}{\begin{Example}\rm}\newcommand{\eex}{\end{Example}}

\newcommand{\pf}{{\bf proof:~}}\newcommand{\proofend}{$\blacksquare$\bigskip}
\newcommand{\bpr}{{\bf proof:~}}\newcommand{\epr}{\hfill$\blacksquare$\\}

\newcommand{\bei}{\begin{itemize}}\newcommand{\eei}{\end{itemize}}
\newcommand{\beq}{\begin{equation}}\newcommand{\eeq}{\end{equation}}
\newcommand{\bem}{\begin{displaymath}}\newcommand{\eem}{\end{displaymath}}
\newcommand{\beqa}{\begin{eqnarray}}\newcommand{\eeqa}{\end{eqnarray}}
\newcommand{\bpm}{\begin{pmatrix}}\newcommand{\epm}{\end{pmatrix}}
\newcommand{\ber}{\begin{array}{l}}\newcommand{\eer}{\end{array}}
\newcommand{\bet}{\begin{tabular}{l}}\newcommand{\eet}{\end{tabular}}

\newcommand{\ra}{\rightarrow}
\newcommand{\di}{\partial}

\newcommand{\li}{~\\ $\bullet$ }

\newcommand{\mC}{\mathbb{C}}
\newcommand{\mP}{\mathbb{P}}
\newcommand{\mZ}{\mathbb{Z}}

\newcommand{\Si}{\Sigma}\newcommand{\tSi}{\tilde{\Sigma}}

\newcommand{\meshn}[7]
{
\put(#1,#2){\vector(1,0){#6}}  \put(#1,#2){\vector(0,1){#7}}
\setcounter{tempx}{#3}  \addtocounter{tempx}{1}  \setcounter{tempy}{#4}  \addtocounter{tempy}{1}
\multiput(#1,#2)(#5,0){\value{tempx}}{\multiput(-1.5,-0.5)(0,#5){\value{tempy}}{.}}
\setcounter{kk}{1}
\setcounter{tempy}{#2} \addtocounter{tempy}{-10}  \setcounter{tempx}{#1} \addtocounter{tempx}{#5}
\multiput(\value{tempx},\value{tempy})(#5,0){#3}{\arabic{kk}\addtocounter{kk}{1}}
\setcounter{kk}{1}
\setcounter{tempy}{#2} \addtocounter{tempy}{#5}  \setcounter{tempx}{#1} \addtocounter{tempx}{-10}
\multiput(\value{tempx},\value{tempy})(0,#5){#4}{\arabic{kk}\addtocounter{kk}{1}}
}
\newcommand{\mesh}[7]
{
\put(#1,#2){\vector(1,0){#6}}  \put(#1,#2){\vector(0,1){#7}}
\setcounter{tempx}{#3}  \addtocounter{tempx}{1}  \setcounter{tempy}{#4}  \addtocounter{tempy}{1}
\multiput(#1,#2)(#5,0){\value{tempx}}{\multiput(-1.5,-0.5)(0,#5){\value{tempy}}{.}}
}
\begin{document}\setcounter{secnumdepth}{6}\setcounter{tocdepth}{2}
\newcounter{temp}\newcounter{kk}\newcounter{iter}\newcounter{Vh}\newcounter{Ch}
\newcounter{tempx}\newcounter{tempy}
\title{Enumeration of singular algebraic curves}
\author{D. Kerner}
\maketitle
\begin{abstract}
We enumerate plane complex algebraic curves of a given degree with one singularity of any given topological type.
Our approach is to compute the homology classes of the corresponding equisingular strata in
the parameter spaces
of plane curves. We suggest a recursive procedure, which is based on the intersection theory combined with liftings
and degenerations. The procedure computes the homology class in question whenever a given singularity type is
defined. Our method does not require the knowledge of all the possible deformations of a given singularity.
\end{abstract}
\tableofcontents
\section{Introduction, definitions and results}
\subsubsection{General settings}
This is an updated and corrected version of \cite{Ker}.\\\\
We work with (complex) algebraic curves of degree $d$ in $\mP^2$.
The parameter space of plane curves is projective space, which we denote by $\mP^D_f$
(here \mbox{$D={d+2\choose{2}}-1$}, the subscript $f$ emphasizes that this is the parameter space of curves defined
by equation $f(x)=0$).
The {\it discriminant}, $\Sigma$, of curves of degree $d$ is the (projective) subvariety of the parameter space,
corresponding to singular curves. There is an equisingular stratification of the discriminant.
An {\it equisingular stratum} is the (quasi-projective) variety of curves with the same topological type of
singularity. The generic point of the discriminant lies in the stratum of nodal curves
(the variety of such curves will be denoted by $\Sigma_{A_1}$). Generic points in the complement,
$\Sigma\backslash\Sigma_{A_1}$, correspond to higher
singularities: the cuspidal ($\Sigma_{A_2}$) or the bi-nodal case ($\Sigma_{A_1,A_1}$).
By further degenerations we get the strata of higher singularities (the precise definition of the equisingular
stratification is in section \ref{SecTypesSingul}):
\beq\scriptstyle
\Sigma=\overline{\Sigma}_{A_1},
~~~\overline{\Sigma}_{A_1}\backslash\Sigma_{A_1}=\overline{\Sigma}_{A_2}\cup\overline{\Sigma}_{A_1,A_1},
~~~\overline{\Sigma}_{A_2}\backslash\Sigma_{A_2}=\overline{\Sigma}_{A_3}\cup\overline{\Sigma}_{A_1,A_2},
~~~\overline{\Sigma}_{A_1,A_1}\backslash\Sigma_{A_1,A_1}=\overline{\Sigma}_{A_1,A_1,A_1}\cup\overline{\Sigma}_{A_1,A_2}
\ldots
\eeq
In this paper we work mainly with isolated singularities. Unless stated otherwise, we mean by $\Sigma_*$ a stratum of
topological type of isolated singularity.

The first natural question is whether a particular stratum is nonempty (i.e. whether a curve of a given degree
can possess prescribed singularities). It still has no complete solution, only lower and upper bounds on the degree
of curves and the codimension of the singularity are known \cite[chapter 5]{Shustin}.

The next question is on the (ir)reducibility, smoothness and dimension of an equisingular stratum. In general a stratum
can be reducible and can have components of different dimensions \cite{GLS,Luengo}. Assuming that an equisingular stratum
is irreducible and has expected dimension we can ask for the degree of its (topological) closure.

Every equisingular stratum is (naturally)
embedded into $\mP^D_f$, so its closure has the homology class (in the integer homology group of the parameter space)
\beq
[\overline{\Sigma}_*]\in H_*(\mP^D_f,\mZ)
\eeq
which is just the needed degree.
Using Poincare duality we obtain the class in cohomology, it will be denoted by the same letter
$[\overline{\Sigma}_*]\in H^*(\mP^D_f,\mZ)$, no confusion should arise.
\subsubsection{The goal of the paper, motivation and main results}
The goal of this paper is to calculate these cohomology classes, for the strata corresponding to curves with one singular
point of a given topological type.
\\
\\
\\
The discriminant, and more generally, varieties of equisingular hypersurfaces, have been a subject of study for a long time.
Already in the 19'th century it was known that the (closure of the) variety of nodal hypersurfaces of degree d in $\mP^n$
is an irreducible projective variety of degree
\beq
(n+1)(d-1)^n
\eeq
Any further progress happens to be difficult.
The present situation is as follows. (This is not a complete/historical review, for a much better description see
\cite{Klei1,Klei2}.)
\li{In 1981 I.Vainsencher\cite{Vain1})} has proven Enriques' formula for the degree of the stratum of cuspidal ($A_2$) curves.
Recently \cite{Vain2} he has calculated the degrees of the strata of multi-nodal hypersurfaces (for the number of nodes $\leq6$).
\li{In 1989 Z.Ran \cite{Ran}} has described an inductive approach to counting plane curves with a prescribed number of nodes.
\li{In 1998 P.Aluffi\cite[section 3]{Alufi}} enumerated cuspidal and bi-nodal hypersurfaces and tacnodal curves.
The results for curves appeared also in \cite[section 2.1]{diFrItz}.
\li{In 1998 L.Caporaso and J.Harris \cite{CapHar1,CapHar2}} have given an algorithm to calculate the degree of the stratum
of curves with any given number of nodes.
\li{S. Kleiman and R. Piene \cite{KleiPien1,KleiPien2}} enumerated multinodal curves with the number of nodes up to 8,
and curves with one triple point and up to 3 nodes.
\li{M.Kazarian in the series of papers \cite{Kaz1,Kaz2,Kaz3,Kaz4}} used topological approach to prove that there
exists a universal formula for the degrees of equisingular strata. It is a (unknown) polynomial in the (relative) Chern classes,
depending on the topological characteristics of the ambient space, singularity type, degree of the hypersurface.
He developed a method for calculations of degrees of strata
for singularities of the given codimension (in particular he presents the answers, for curves, for equisingular strata up to
co-dimension 7).
The drawback of his method is that before one starts to calculate, one should classify all the singularities of the given
codimension. So, for example, it is extremely difficult to calculate the degree of $[\Sigma_{A_{k}}]$, for a given large $k$
(meanwhile there is no corresponding classification). Another omission is that it is unclear
what is the "domain of universality" of the formulae (for small degrees, relatively to the codimension
of the singularity, the stratum is reducible and the universality does not hold).
\li{} Another proof of the existence of the universal formula was given in 2004 by A.K.Liu \cite{Liu}
(by algebraic methods). He, however, does not give any concrete methods of calculation.
\li One should also mention numerous strong results by  Feh\'er, Rim\'anyi,
Lascoux, Fulton and many others.

Our motivation was to solve the enumeration problem by use of the classical intersection theory approach and to avoid the
difficulties that occur in other approaches. Our main result is {\it the algorithm, which enumerates algebraic curves
(gives the degree of the corresponding equisingular stratum) with one singular
point of any (given) topological type, provided we know the normal form of the type}.

 The advantage of our method is that it can be directly
implemented for every particular class of singularity, without any preliminary classification. The only initial data to
start is the normal form.
In fact, during the degeneration procedure (described in details in section \ref{SecTheDegen}), we only need to know
some singularity types
adjacent to the given singularity in
codimension 1. These types appear explicitly in course of calculations.

In case of linear-singularities (defined in section \ref{SecDefinitPropertLinearStrata}) the method gives immediate answer, in other cases it provides an algorithm
(which is quite efficient and can be programmed).
We discuss the case of Newton-non-degenerate singularities in full details. For other cases the method generalizes
naturally, however calculations become tedious and do not bring anything new from the conceptual point of view. So,
instead of treating the Newton-degenerate case in full generality, we consider some typical classes of such
singularities, where all the ingredients of the generalization are demonstrated.
\subsubsection{Description of the method}
The main idea of the method is to work with complete (or at least locally complete) intersections of hypersurfaces.
For this we lift the given equisingular stratum (which naturally lies in $\mP^D_f$, the parameter space of plane curves)
to a bigger space. For example for unibranched singularity we define:
\beq
\tilde\Sigma=\overline{\left\{(x,l,f)|~\ber \mbox{The~curve~defined~by~}\Big(f(x)=0\Big)~\mbox{has~the~prescribed}\\
\mbox{singularity~at~the point~}x\mbox{~with~the tangent~line}~l\eer\right\}}\subset
\mP^2_x\times(\mP^2_l)^*\times\mP^D_f
\eeq
Here $\mP^2_x$ is the plane, the subscript $x$ emphasizes that the point of the plane is denoted by $x$,
the same for $(\mP^2_l)^*$. To avoid messy notations from now on we mean by $\Sigma$ (or $\tilde\Sigma$) the closure
of a (lifted) equisingular stratum

The cohomology class of the lifted version of a stratum is often easier to calculate (e.g. for ordinary multiple point the lifted
stratum is a globally complete intersection).
Now, the (co)homology class is not just a number, but a polynomial (in the generators of the cohomology ring
of the bigger embedding space), and hence we have a multidegree: $[\tilde\Sigma]$.
This of course provides
much more information about a particular stratum. We present some of the calculated multi-degrees in the table below.
Once the class $[\tilde\Sigma]$ is calculated, the cohomology class of the
(closure of) original stratum is obtained as follows. The projection $\tilde\Sigma\stackrel{\pi}{\rightarrow}\Sigma$ (generically
1:1) induces a projection of homology. Therefore (by Poincare duality) we have a map in cohomology (the Gysin
homomorphism), which gives the (cohomology) class of $\Sigma$.

Even if the lifted version is not a globally complete intersection we could hope to define it by some standard
(well studied)
conditions (e.g. conditions of proportionality of two tensors, matrices of the given rank \cite{Ful}).
In this way we could express the cohomology class of the stratum in terms of the standard classes.
This is the case for linear singularities\footnote{e.g. $A_{k\le3}~D_{k\le6},E_{k\le8},~X_9,~J_{10},~Z_{k\le13}$ etc.}
 (defined in section \ref{SecDefinitPropertLinearStrata}).
In section \ref{SecLinearStrata} we prove the following theorem:

{\it Every linear stratum $\Sigma$ can be lifted to a variety $\tilde\Sigma$, which is a locally complete intersection.
$\tilde\Sigma$
is a variety of specific type, defined by transversal intersection of conditions, each of them being proportionality of symmetric
multi-forms: $\Omega^{(p_i)}_i\sim \omega^{(p_i)}_i$. Correspondingly, the cohomology class is calculated
as a product of these conditions.}

In many cases, even the class of the lifted version cannot be easily calculated (the case of non-linear-strata,
section \ref{SecImpossibilityOfTheDirectApproach}).
Instead we use the method of degeneration. Namely, we add some restrictions to the defining conditions of the
singularity, to obtain a singularity of higher type that has linear stratum and so is simple to work with. Or, geometrically,
we intersect
the (lifted) stratum with a divisor (or a cycle) so that the cohomology class of the intersection is easy to calculate.
And from the cohomology class of the intersection we can get some information about the cohomology class of the
original stratum. Actually, in all cases considered in the paper, the class of the original stratum is recovered
uniquely (we discuss the invertibility of degeneration in section \ref{ParDegenCondit}).

More formally, for every singularity type we start from, we build the degeneration tree with the following properties:
\li The tree is finite.
The root of the tree is the initial stratum, while the leaves correspond
to some specific linear strata for which the enumeration problem is solved explicitly (as described above).
Other vertices correspond to some intermediate singularity types.
\li To each vertex there are assigned (oriented) edges corresponding to a minimal (codimension 1) degeneration
of the current singularity.
The neighbors of a vertex (in the positive direction) correspond to the intermediate
strata that appear as the result of degeneration. The edges are decorated by the multiplicities with
which the corresponding strata appear in the degeneration.
\li Each vertex is a root of a (finite) subtree. This subtree (its edges and vertices) is defined by
its root only (i.e. by its singularity type).
\li Each vertex (together with its outgoing edges) gives a
linear equation in the cohomology rings of the embedding space, relating the cohomology class of the initial
stratum to those of its neighbors. This equation  always has a unique solution.

Once the degenerating tree is constructed, the enumerative problem for the root of the tree (the initial singularity)
is reduced to a collection of much simpler problems for the leaves of the tree (linear singularities).

In this procedure several issues constantly repeat themselves e.g. transversality of the defining conditions,
multiplicity of the obtained variety, residual pieces "at infinity". We describe them in details in section
\ref{SubSubSecNonLinearStrata}.

As in every problem whose final answer is an explicit numerical formula, it is better to have several different solutions
leading to the same result (in addition to the rigorous proofs). Many of the results in the paper are derived in several
independent ways (e.g. the cases of $A_{k<4},~D_{k<7}$). Another kind of check is provided by the known results for the
cases of small degrees. This happens in particular when the curve is reducible
(a curve of degree $d$ with an $A_{(d-1)(d-2)+1}$ singularity necessarily decomposes into two maximally tangent curves).
 The degree of the strata in these cases can be calculated directly.
And (as it should be) we obtain the same result.
In some cases there are additional consistency checks (e.g. the multidegree should be symmetric with respect to some
variables or should not contain some particular monomials). They help to avoid lengthy calculations of some unknown
parameters.

We should note, that, though our approach works well in every particular case (e.g. for $A_k$ with every particular $k$),
it
does not allow to obtain general (universal) formulae (e.g. the general formula for all the $A_k$'s).
\subsubsection{The range of the universality}
An important question is the applicability range of the obtained formulae (or of the algorithm). For example
 a cubic cannot posses
 the (isolated) $A_k$ singularity for $k\geq4$. Correspondingly, the answers in the table are non-positive for
 $d=3,~k\geq4$. Even worse, in many cases of low degree the curve can possess the prescribed (isolated) singularity
 but (due to the low degree) various
coincidences/degenerations can happen. The answer in this case is not described by the universal formulae, in particular the algorithm
gives wrong (often negative) degrees\footnote{One example of this is a quartic with $A_7$}.
As is stated in \cite{Kaz1}-\cite{Kaz4}, the universality exists when the degree of the hypersurface is sufficiently high.
Our algorithm provides explicit sufficient bounds on the degree, one of them is:
\beq
\text{order of determinacy}\le d+1
\eeq
Another sufficient bound seems to exist: {\it codimension~of~singularity$<$2d-1}. We discuss it in section
\ref{SecNonUniversalDomain}.

The non-universality domain must be treated separately. Due to its non-universal nature, there probably does not exist any
general formula (even for a given type of singularity). Rather, one should calculate case-by-case (i.e. for a given
degree of curve and a given singularity type). This problem is solved by degeneration techniques. We illustrate
this by the example of quartics with $A_7$ in section \ref{SecCalcInNonUnivDom}.
\subsubsection{Degrees of the equisingular strata}
In this paper we widely use notations for singularities, taken from \cite[chapter II]{AGV}. The classes of singularities are denoted
by letters or by their normal form or by the jets of their normal form, e.g.
 $A_k:~~jet_k(f)=x^2_2,~~X_9:~~jet_3(f)=0$.
\\
The algorithm calculates the degrees of the strata, or, alternatively, the cohomology classes
$[\Sigma]\in H^{2\rm{codim}}(\mP^D_f)$. Here codim is the codimension of a given stratum in the parameter space.
For simple singularities
(A-D-E) it equals the Milnor number $\mu$, in general codim$=\mu-$modality.

As the simplest example of the use of the method we present the following results:
\bec
The degrees of the equisingular strata for some low codimension singularity types are given in the tables below
\eec
\begin{tabular}{|@{\small}>{$}l<{$}@{\small}|@{\small}>{~$}l<{$}@{\small}|}\hline
A_1&3(d-1)^2\\\hline
A_2&12(d-1)(d-2)\\\hline
A_3&2(25d^2-96d+84)\\\hline
A_4&60(d-3)(3d-5)\\\hline
A_5&18(35d^2-190d+239)\\\hline
A_6&7(316d^2-1935d+2769)\\\hline
A_7&12(651d^2-4400d+7002)\\\hline
\end{tabular}
\begin{tabular}{|@{\small}>{$}l<{$}@{\small}|@{\small}>{~$}l<{$}@{\small}|}\hline
D_4&15(d-2)^2\\\hline
D_5&12(d-2)(7d-19)\\\hline
D_6&14(16d^2-87d+114)\\\hline
D_7&48(15d^2-92d+135)\\\hline
D_8&2(2025d^2-12438d+18392)\\\hline
\end{tabular}
\begin{tabular}{|@{\small}>{$}l<{$}@{\small}|@{\small}>{~$}l<{$}@{\small}|}\hline
W_{12}&24(d-4)(13d-44)\\\hline
W_{13}&3(4128-2132d+273d^2)\\\hline
W_{17}&2584d^2-23664d+53259\\\hline
W_{18}&3(1173d^2-11288d+26608)\\\hline
W_{24}&144(462d^2-5173d+14135)\\\hline
W_{25}&24(6325d^2-72732d+204362)\\\hline
W_{1,0}&12(105d^2-868d+1773)\\\hline
\end{tabular}
\\
\begin{tabular}{|@{\small}>{$}l<{$}@{\small}|@{\small}>{~$}l<{$}@{\small}|}\hline
J_{10}&2(385d^2-2715d+4611)\\\hline
Z_{11}&18(267 - 154d + 22d^2)\\\hline
Z_{13}&3(286d^2-2148d+3977)\\\hline
X_{1,1}&30(11d-39)(d-3)\\\hline
X_{1,2}&3(242d-1727d+3035)\\\hline
\end{tabular}
\begin{tabular}{|@{\small}>{$}l<{$}@{\small}|@{\small}>{~$}l<{$}@{\small}|}\hline
E_6&21(d-3)(4d-9)\\\hline
E_7&252d^2-1464d+2079\\\hline
E_8&9(45d^2-288d+448)\\\hline
E_{12}&12(442d^2-3540d+6851)\\\hline
E_{13}&8(1820d^2-15639d+32353)\\\hline
\end{tabular}

Many more formulae are in given Appendix.
\subsubsection{Multi-degrees of the equisingular strata}
In the table below we present the multi-degrees (cohomology classes) of the lifted equisingular strata,
for the lifting:
\beq\scriptstyle
\tilde\Sigma=\left\{(x,\{l_i\}_i,f)|~\ber \mbox{The~curve}~\{f(x)=0\}~\mbox{has~the~prescribed singularity}\\\mbox{at~the point x~with~lines}~\{l_i\}_i\mbox{ tangent to the branches}\eer\right\}
\subset\mP^2_x\times\Big\{(\mP^2_{l_i})^*\Big\}_i\times\mP^D_f
\eeq
For unibranched singularities the projection $\tilde\Sigma\mapsto\Sigma$ is generically 1:1; correspondingly, to find the actual degree of
the stratum, one must extract the coefficient of the generator of $H^4(\mP^2_x)\otimes H^4((\mP^2_l)^*)$ (i.e. the coefficient of $x^2l^2$).

Most of the low codimension singularities are either uni-branched or with one tangent line common to several branches.
The simplest exceptions are e.g. ordinary multiple points ($A_1,D_4,X_9,\dots$), here there is no unique tangent line.
In this case we define the lifted variety as:
\beq
\tilde\Sigma=\left\{(x,f)|~f^{(k)}(x)=0\right\}
\eeq
As in the rest of the paper, $X$ denotes the generator of the cohomology ring of $\mP^2_x$.
Similarly: $H^*((\mP^2_l)^*)=\mZ[L]/(L^3),~~H^*(\mP^D_f)=\mZ[F]/(F^{D+1})$.
\bec
The cohomology classes (multi-degrees) of the equisingular strata for some singularity types are given in the tables below
\eec
\begin{tabular}{|@{}>{\footnotesize}l@{}|>{\footnotesize}l@{}|@{}>{\footnotesize}l@{}|}\hline
Ordinary point  of multiplicity $p$ (e.g. $A_1:p=1, D_4:p=2,X_9:p=3$) & $Q^{p+2\choose{2}}$ & $Q=(d-p)X+F$\\\hline
\end{tabular}
\\
\hspace{0cm}{
\begin{tabular}{|@{}>{\scriptsize$}l<{$}@{}|>{\scriptsize$}l<{$}@{}|@{}>{\scriptsize$}l<{$}@{}|}\hline
A_2 & (L+X)(Q^5+2lQ^4+4l^2Q^3) & Q=(d-2)X+F \\\hline
A_3 & (L+X)(Q^6+(5L+2X)Q^5+(16l^2+11lX+X^2)Q^4+42lX^2Q^3) & Q=(d-3)X+F\\\hline
A_4 & 2(L+X)(Q^7+6lQ^6+24l^2Q^5) & Q=(d-3)X+F\\\hline
A_5 & 2(L+X)(2Q^8+(15L+10X)Q^7+(135l^2+83X^2)Q^6+351lX^2Q^5) & Q=(d-4)X+F\\\hline
A_6 & (L+X)(9Q^9+(79L+106X)Q^8+(1225l^2+1328X^2)Q^7+6958lX^2Q^6) & Q=(d-5)X+F\\\hline
A_7 & (L+X)(21Q^{10}+(217l+208X)Q^9+3(3165l^2+2658X^2)Q^8+15324lX^2Q^7) & Q=(d-5)X+F\\\hline
\end{tabular}}
\\
\hspace{-0cm}{
\begin{tabular}{|@{}>{\scriptsize$}l<{$}@{}|>{\scriptsize$}l<{$}@{}|@{}>{\scriptsize$}l<{$}@{}|}\hline
D_5 & (L+X)(Q^8+(4L+X)Q^7+(12l^2+6lX)Q^6+24Xl^2Q^5) & Q=(d-3)X+F\\\hline
D_6 & (L+X)(Q^9+(8L+5X)Q^{8}+(38l^2+44lX+9X^2)Q^7+(224Xl^2+84X^2l)) & Q=(d-4)X+F\\\hline
D_7 & 2(L+X)(Q^{10}+2(5L+6X)Q^9+3(59X^2+58l^2)Q^8+1200Xl^2Q^7) & Q=(d-5)X+F\\\hline
D_8 & 2(L+X)(3Q^{11}+45LQ^{10}+(418l^2+74X^2)Q^9+1040Xl^2Q^8) & Q=(d-4)X+F\\\hline
\end{tabular}}
\\
\begin{tabular}{|@{}>{\scriptsize$}l<{$}@{}|>{\scriptsize$}l<{$}@{}|@{}>{\scriptsize$}l<{$}@{}|}\hline
E_6 & (L+X)(Q^9+3lQ^8+9l^2Q^7) & Q=(d-3)X+F\\\hline
E_7 & (L+X)\Big(\ber Q^{10}+(31l^2+38lX+10X^2)Q^8+\\(7L+5X)Q^9+(173l^2X+82lX^2)Q^7\eer\Big) & Q=(d-4)X+F\\\hline
E_8 & (L+X)\Big(\ber Q^{11}+(51l^2+21lX+X^2)Q^9+\\(9L+2X)Q^{10}+(132l^2X+12lX^2)Q^8\eer\Big) & Q=(d-4)X+F\\\hline
E_{12} &
4(L+X)(Q^{14}+(17L+3X)Q^{13}+(220l^2+57X^2)Q^{12}+603lX^2Q^{11}) &
Q=(d-5)X+F\\\hline E_{13} &
4(L+X)(2Q^{15}+(40l-3X)Q^{14}+2(197l^2-30X^2)Q^{13}-684lX^2Q^{12})
& Q=(d-5)X+F\\\hline
J_{10} & (L+X)\Big(\ber Q^{12}+(111l^2+116lX+26X^2)Q^{10}+\\(14L+8X)Q^{11}+(936l^2X+386lX^2)Q^9\eer\Big)
& Q=(d-5)X+F\\\hline J_{2,1} &
2(L+X)\Big(\ber 2Q^{13}+(31L+37X)Q^{12}+\\6(137l^2+143X^2)Q^{11}+8931lX^2Q^{10}\eer\Big)
& Q=(d-6)X+F\\\hline
\end{tabular}
\\
\begin{tabular}{|@{}>{\scriptsize$}l<{$}@{}|>{\scriptsize$}l<{$}@{}|@{}>{\scriptsize$}l<{$}@{}|}\hline
W_{24} &
8(L+X)(Q^{23}+(36L+11X)Q^{22}+10(111l^2+46X^2)Q^{21}+10116lX^2Q^{20})
& Q=(d-7)X+F\\\hline W_{25} &
8(L+X)(2Q^{24}+(75L+16X)Q^{23}+10(2157l^2+682X^2)Q^{22}+15489lX^2Q^{21})
& Q=(d-7)X+F\\\hline
\end{tabular}
\subsubsection{Organization of material}
In section 2 we fix the notations, recall some important notions and describe the method. Then we
demonstrate the method in the simplest cases: complete intersections (ordinary multiple point),
liftings of linear singularities (the cuspidal case), degenerations (the cuspidal case).

In section 3 we consider the case of linear singularities. We first introduce linear
singularities and explore their properties. Then we solve the problem for linear strata
by lifting them to bigger spaces
(the result is formulated in Lemma \ref{LemmaLinear}).
The general formula for the degree (though it exists), must be very complicated, therefore we treat each case separately.
Several explicit calculations are in Appendix \ref{SecExplicitCalculationLinearStrata}.

At the end of the section 3 we demonstrate another method: degeneration. In case of linear strata it is just another
check of the already obtained answers, in the case of non-linear-strata it is the only currently working method.

In section 4 we deal with non-linear-strata (using degeneration method).

We describe the algorithm of degeneration and use it to calculate the degrees in several cases.
The calculations in this case are much more involved, we solve several cases in Appendix \ref{SecExplicitCalculationNonLinearStrata}.

In sections 3 and 4 we consider Newton-non-degenerate singularities (the definition is in section \ref{SecTypesSingul}).
In subsection 4.4 we consider Newton-degenerate
singularities, discuss the generalization of the method and solve some typical cases.

In section 5 we discuss the range of universality and explain how to solve the problem in the non-universal
domain (when the degree of the curve is low compared to
the codimension of the stratum). We demonstrate how the method works by the example of quartics with $A_7$.

In the body of the paper we frequently use several technical results on cohomology classes of particular conditions
(proportionality of tensors, tangency of curves, specific degenerations etc.). These
conditions are described in Appendix, where all the cohomology classes are explicitly obtained.
\subsubsection{Acknowledgements}
I wish to thank my supervisor E. Shustin for his constant help and support during the work. The solution of the problem
would have been impossible without his geometric ideas and great patience, when answering my endless technical questions.

The work was initiated during the author's stay at the Universit\"at Kaiserslautern.
Significant part of this work was done during the seminar on "Computational Algebraic Geometry" at the
"Mathematische Forschungsinstitut-Oberwolfach". I would like to thank these institutions for excellent working conditions.

I would like to thank R.Piene for inspiring discussion and numerous remarks, which significantly improved the paper.
I would like to thank the referee for pointing out various misprints and omissions of the first version of manuscript.
The research was constantly supported by the Hermann-Minkowski Minerva Center for Geometry at Tel Aviv
University and by the Israel Science Foundation grant, no: 465/04.

Finally I would like to thank M.Kazarian for his valuable comments. They helped to clarify various fine aspects
and to correct numerical errors that occurred in the previous versions of the paper.
\section{Preliminaries}
\subsection{Some definitions}
\subsubsection{On variables}\label{SectVar}
In this paper we deal with many (co)homology classes of various varieties, embedded into various (products of)
projective spaces.
To simplify the formulae we adopt the following notation. If we denote the point in the space
$\mP^n_x$ by the letter $x$ then the homogeneous coordinates are $(x_0,\dots,x_n)$. The generator of the cohomology
ring of this $\mP^n_x$ is denoted by the upper-case letter $X$, so
that $H^*(\mP^n_x)=\mZ[X]/(X^{n+1})$. By the same letter we also
denote the hyperplane class in homology of $\mP^n_x$. Since it is always clear,
where we speak about coordinates and where about (co)homology classes no confusion arises.

To demonstrate this, consider the hypersurface:
\beq
V=\{(x,y,f)|~f(x,y)=0\}\subset\mP^n_{x}\times\mP^n_{y}\times\mP^D_f
\eeq
Here $f$ is a bi-homogeneous polynomial of bi-degree $d_x,d_y$ in homogeneous coordinates \mbox{$(x_0,\dots,x_n)$},
$(y_0,\dots,y_n)$,
the coefficients of $f$ are the homogeneous coordinates in the parameter space $\mP^D_f$.
The cohomology class of this hypersurface is:
\beq\label{DemostratCohomClasses}
[V]=d_xX+d_yY+F\in H^2(\mP^n_{x}\times\mP^n_{y}\times\mP^D_f)
\eeq
A (projective)
line through the point $x$ is defined by a 1-form $l$ (so that: $l\in(\mP^2_l)^*,~l(x)=0)$. Correspondingly
the generator of $H^*((\mP^2_l)^*)$ is denoted by $L$.

We will often work with symmetric $p-$forms $\Omega^{p}:S^p V\rightarrow \mC$
(here $V$ is an $3$ dimensional vector space).
Thinking of the form as of a symmetric tensor with $p$ indices ($\Omega^{(p)}_{i_1,\dots,i_p}$), we often write
$\Omega^{(p)}(\underbrace{x,\dots,x}_{k})$ as a shorthand for the tensor, multiplied $k$ times by the point $x\in\mP^2_x$
(which is considered here as a vector in $V$):
\beq
\Omega^{(p)}(\underbrace{x,\dots,x}_{k}):=\sum_{0\le i_1,\dots,i_k\le2}\Omega^{(p)}_{i_1,\dots,i_p}x_{i_1}\dots x_{i_k}
\eeq
So, for example, the expression $\Omega^{(p)}(x)$ is a $(p-1)-$form. Unless stated otherwise, we assume the
symmetric form $\Omega^{(p)}$ to be generic (in particular non-degenerate, i.e. the corresponding hypersurface
$\{\Omega^{(p)}(\underbrace{x,\dots,x}_{p})=0\}\subset\mP_x^n$ is smooth).

Symmetric forms will typically occur as tensors of derivatives of order $p$: $f^{(p)}$. Sometimes, to emphasize
the point at which the derivatives are calculated we will assign it. So, e.g.
$f^{(p)}|_x(\underbrace{y,\dots,y}_k)$ means: the tensor of derivatives of the
$p$'th order, calculated at the point $x$, and contracted $k$ times with $y$. Usually it will be evident, at which point
the derivative is calculated, in such cases we will omit the subscript $x$.

Thinking of a symmetric $p-$form as of a tensor in $p$ indices, one can consider the wedge product of the $p-$form
and a one-form:
\beq
(\Omega^{(p)}\wedge l)_{i_1\dots i_p j}:=\Omega^{(p)}_{i_1\dots i_p}l_j-\Omega^{(p)}_{i_1\dots i_{p-1}j}l_{i_p}
\eeq
The products of this type will constantly occur in enumeration of linear singularities. Note that the wedge product of a
one-form with itself is trivially zero: $l\wedge l\equiv0$. The inverse statement is also true:
\bec
let $\Omega^{(p)}$ be a symmetric $p-$form satisfying $\Omega^{(p)}\wedge l=0$. Then there exists a symmetric
$(p-1)$ form $\tilde\Omega^{(p-1)}$ such that: $\Omega^{(p)}=SYM(\tilde\Omega^{(p-1)},l)$
\eec
Here $SYM$ means complete symmetrization of indices.
This is a standard result we omit the proof.

Throughout the paper we will tacitly assume Euler identity for a homogeneous
polynomial of degree $d$
\beq
\sum_ix_i\di_if=df
\eeq
and its consequences (e.g. $\sum_ix_i\di_i\di_jf=(d-1)\di_jf$). So, for example the nodal point
(defined by $f^{(1)}|_x=0$) can be defined by $f^{(p)}|_x(\underbrace{x,\dots,x}_{p-1})=0$.
\subsubsection{On the singularities}\label{SecTypesSingul}
For completeness we recall some notions related to singularities of curves \cite{Shustin,AGV}.
\bed
Let $(f,x)\subset(\mC^2_x,x)$ and $(g,y)\subset(\mC^2_y,y)$ be two germs of isolated curve singularities. They are
 topologically (analytically) equivalent if there exist a homeomorphism (local analytic map)
 $(\mC^2_x,x)\mapsto(\mC^2_y,y)$
 mapping $(f,x)$ to $(g,y)$. The corresponding
equivalence class is called topological (analytical) singularity type. The variety of points (in the parameter space),
corresponding to curves with singularity of the same (topological/analytical) type is called the {\bf equisingular stratum}
\eed
\parbox{12.5cm}
{The topological type of a singularity of a plane curve is (uniquely) defined by the minimal embedded resolution tree
and the multiplicities of the infinitely near points \cite[section 4.3]{Shustin}. The singularity is resolved by a
series of blow-ups. At each step we record the points of intersection of the strict transform with the
exceptional divisor. To each such 'infinitely near point' we assign its multiplicity. For purposes of our paper we
assign to each infinitely
near point the degree of intersection of the strict transform with the exceptional divisor of the corresponding blowup.}
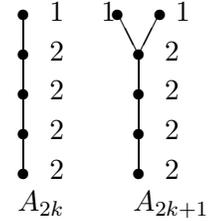
\begin{picture}(0,0)(-20,-20)
\put(0,-30){\line(0,1){60}}  \multiput(0,-30)(0,15){5}{\circle*{4}}  \multiput(10,-32)(0,15){4}{2} \put(10,28){1}
\put(-2,-44){$A_{2k}$}
\end{picture}
\begin{picture}(0,0)(-60,-20)
\put(0,-30){\line(0,1){45}}  \multiput(0,-30)(0,15){4}{\circle*{4}}  \multiput(10,-32)(0,15){4}{2}
\put(-2,-44){$A_{2k+1}$}
\put(0,15){\line(1,2){7}}   \put(0,15){\line(-1,2){7}}
\put(-8,30){\circle*{4}} \put(-14,28){1}  \put(8,30){\circle*{4}}   \put(14,28){1}
\end{picture}
\\
\\
For each topological singularity type there are many particular representatives. Among them one chooses (a generic
and most simple) one which is called {\it the normal form} for this type. For example for several simplest types
(all the notations are from \cite{AGV}, we ignore the moduli):
\beq\scriptstyle\ber
A_k:x^2_2+x^{k+1}_1,~~D_k:x^2_2x_1+x^{k-1}_1,~~E_{6k}:x^3_2+x^{3k+1}_1,~~E_{6k+1}:x^3_2+x_2x^{2k+1}_1,~~
E_{6k+2}:x^3_2+x^{3k+2}_1\\
J_{k\ge1,i\ge0}:x^3_2+x^2_2x^k_1+x^{3k+i}_1,~~Z_{6k-1}:x^3_2x_1+x^{3k-1}_1,~~
Z_{6k}:x^3_2x_1+x_2x^{2k}_1,~~Z_{6k+1}:x^3_2x_1+x^{3k}_1\\
X_{k\ge1,i\ge0}:x^4_2+x^3_2x^k_1+x^2_2x^{2k}_1+x^{4k+i}_1,~~W_{12k}:x^4_2+x^{4k+1}_1,~~W_{12k+1}:x^4_2+x_2x^{3k+1}_1
\eer\eeq
In the case of simple singularities (A,D,E) the curve can be brought to the normal form by locally analytic transformations.
Correspondingly, the topological and analytical approaches coincide.
However, for higher singularities this is not the case\footnote{For example, any two ordinary k-tuple points ($k\geq4$) are topologically equivalent,
but there is a continuum of different analytical types (e.g. for $k=4$ the cross-ratio of tangent lines is fixed by analytic
transformations) },
various moduli appear. From the topological point of view, singularities with different moduli belong to the same class.
Therefore in
this paper the moduli are always assumed generic and ignored. In particular we usually assume that the nonzero coefficients of monomials
in the normal form are
equal to 1.

Using the normal form one can draw
the Newton diagram of the singularity. Namely, one marks the points corresponding to non-vanishing monomials,
and takes the convex hull.

From the Newton diagram of the singularity, one immediately records monomials that should be absent in the
normal form. We formalize this fact in the notion of the {\it preliminary form} (for a given topological type and
a given normal form), which will be very useful throughout the paper:
\bed
The germ of singular curve is in a preliminary form if it has the given singularity type and its Newton diagram
coincides with the Newton diagram of the normal form of this singularity.
\eed
For enumeration of Newton-non-degenerate singularities we need the preliminary form only.

Note, that the curve that was brought to a preliminary form is not necessarily in the normal form.
The simplest example is just the ordinary multiple point: here the preliminary form means that the singularity
is at the origin and the local germ is $f=\Omega^{(p)}$ (a generic homogeneous $p-$form). The normal form of
the ordinary multiple point is more special: $f=x^p_1+x^p_2$.

The distinction is even more serious: as the Newton diagram does not always specifies the singularity type,
different types could have the same preliminary form. The singularities, for which the Newton diagram
specifies the type completely are called Newton-non-degenerate.
\bed
The singularity is (generalised) Newton-non-degenerate if the truncation of
the corresponding polynomial to every segment of the Newton polygon is non-degenerate\footnote{the truncated polynomial has no
singular points in the torus $(C^*)^2$} or the singularity can be brought to such a form by a
locally analytic transformation.
\eed

To each equisingular stratum there are adjacent strata of higher singularities.
For example: $\Sigma_{D_{k+1}}\subset\overline\Sigma_{A_k}$.
We will primarily be
interested in the codimension one adjacency, i.e. the strata that can be reached by just one degeneration.
Many tables of such adjacencies are given in \cite{AGV}. Recently the question of, whether the two strata with
the given Enriques diagrams are adjacent, was solved \cite{AR}. In each particular case this question can be
answered by the careful analysis of Newton diagram of the singularity.

There are examples when adjacency depends on moduli \cite{Pham}. However, if one chooses moduli generically
(and we always
do that) then adjacency is completely fixed by the topological type.

We denote (the closure of) an equisingular stratum by $\Sigma$ (in each case it
will be clear which singularity is meant).
\subsection{Main issues of the method}
 Our goal is to calculate the (co)homology class of the equisingular strata.
We explain here the main issues of the method. All the steps are illustrated
by simple examples in subsection \ref{SecExamples}. Consulting them will greatly simplify the reading.
\subsubsection{Complete, locally complete and not locally complete intersections}\label{SecComplLocComplNotLocCompl}
As stated in the introduction we use classical intersection theory. We try to define the varieties by some
explicit equations (i.e. to represent them as the intersection of hypersurfaces). Then we try to relate
the cohomology class of the product to that of the needed variety. Several cases occur:
\li Globally complete intersection. The hypersurfaces intersect transversally, therefore the cohomology class
is just the product of classes of the hypersurfaces (the later being fixed by the (multi)degrees).
This happens in the case of ordinary multiple point (\ref{SecOrdMultPoint}).
\li Locally complete intersections are locally defined by a transversal set of hypersurfaces, but not globally.
To calculate the cohomology classes of such varieties one chooses an open set where the variety is
a complete intersection. Choose a set of hypersurfaces $\{V_i\}_{i=1}^k$ that defines the variety locally.
Now consider the (globally defined) variety: $\cap_{i=1}^k V_i$. It coincides with the original variety in an open set,
however out outside this set ("at infinity") some residual varieties can appear.
 Then the cohomology class of the original variety is calculated as:
\beq
\Pi_{i=1}^k[V_i]-[\mbox{the residual pieces "at infinity"}]
\eeq
The residual pieces at infinity usually come with some multiplicities, these are fixed in each particular case by
the analysis of the equations (see the simplest example of nodal stratum in subsection \ref{SubSubSectionNodal}).
\li For not locally complete intersections, the above method does not work. Instead we can try
to lift them to a bigger space (to make them locally complete intersection), or to degenerate them to some (simple)
varieties of higher codimension.
\subsubsection{Liftings}\label{SecLiftings}
Most equisingular strata are not globally complete intersections.
To calculate their cohomology classes, the first thing we do is to lift the strata to a bigger space. Namely, instead of
a particular point of the stratum: $f\in\mP_f^D$ (that defines a plane curve: $\{f(x)=0\}\subset\mP_x^2$), we
consider the
triples: $(x,l,f)$. Here, $x\in\mP^2_x$ is the singular point of the curve, $l\in(\mP^2_l)^*$ is a one-form that
defines a line through $x$, tangent to one of the branches of the singularity. In this way we define the
{\it lifted stratum}
\beq\scriptstyle
\tilde\Sigma(x,l)=\overline{\{(f,x,l)|\ber\mbox{the curve~defined~by~}(f(x)=0)~\mbox{has}\\\mbox{prescribed~singularity~at~}x,
~\mbox{with~tangent~line~}l\eer\}}
\subset\mP^D_f\times\mP_x^2\times(\mP_l^2)^*
\eeq
This lifting suites in a uni-branched case or when the different branches have the same tangent line. In general
one should of course lift further to take into account all the tangents.

The projection $\mP^D_f\times\mP_x^2\times(\mP_l^2)^*\stackrel{\pi}{\rightarrow}\mP^D_f$, restricted to $\tilde\Sigma$
is generically 1:1. We have the induced map in homology:
$H_k(\mP^D_f\times\mP_x^2\times(\mP_l^2)^*)\stackrel{\pi_*}{\rightarrow}H_k(\mP^D_f)$ and (by Poincare duality) the
 map in cohomology: $H^{k+4}(\mP^D_f\times\mP_x^2\times(\mP_l^2)^*)\stackrel{\pi_*}{\rightarrow}H^k(\mP^D_f)$ (the
 Gysin homomorphism).

 We will need the precise way, the Gysin homomorphism acts. Start from the projection in homology. Let $V$ be a linear
 subspace in $\mP^D_f\times\mP_x^2\times(\mP_l^2)^*$ (a representative of a homology class). Then (by dimensionality)
 $\pi_*[V]=0$,
 unless the projection: $V\stackrel{\pi}{\rightarrow}\pi(V)$ is 1:1 (in which case it is an isomorphism).
Thus $\pi_*$ acts as identity on: $H_k(\mP^D_f)\otimes H_0(\mP_x^2\times(\mP_l^2)^*)$ and as zero on other homology
groups.
Correspondingly, it acts as identity on $H^k(\mP^D_f)\otimes H^4(\mP_x^2\times(\mP_l^2)^*)$ and sends all other
cohomology classes to zero.

From the calculational point of view, to obtain the cohomology class of $\Sigma$ we just need to extract the
coefficient of $X^2L^2$ from the cohomology class of $\tilde\Sigma(x,l)$ (here $X,L$ are the generators of the
cohomology rings of $\mP_x^2$ and $(\mP_l^2)^*$.

In general, one has to lift further, to take into account other parameters of singularity. Then to go
back Gysin projection is applied again.

It may happen that the projection $\tSi\ra\Si$ is a (finite) covering or a fibration with projective spaces
as the fiber.

Summarizing: the cohomology class of a stratum $\Sigma$ is completely fixed by the class of its lifted version
$\tilde\Sigma$. Therefore, from now on we will be interested in $[\tilde\Sigma]$ only.
The simplest case of lifting (for cuspidal curves) is considered in \ref{SecCusp}.
\subsubsection{Linear strata (Newton non-degenerate singularities)}
The cohomology class of $\tilde\Sigma$ is often easier to calculate. In particular for a big class of
(Newton-non-degenerate) singularities, $\tilde\Sigma$ is a locally complete intersection, and we can write explicit,
global, covariant equations defining the lifted strata.
The strata/singularities that can be enumerated in such a simple way are called {\it linear}
(the precise definition is in section \ref{SecDefinitPropertLinearStrata}).

There are two reasons for this simplicity. The first is a simple algebraic fact: every homogeneous form of rank two
fully factorizes (in other words, any homogeneous polynomial in two variables is totally reducible).
The second: for linear singularities the defining equations of the corresponding lifted strata are linear in functions or
their derivatives.
The simplest (nontrivial) linear singularity is the cusp ($A_2$), its enumeration is done in subsection \ref{SecCusp}.

As stated in Lemma \ref{LemmaLinear}, the lifted strata can be represented by mutually transversal conditions,
each condition amounts to proportionality of two (symmetric) tensors. The cohomology class of this last condition
is classically known, it is given in Appendix \ref{TensorProportionality}. In total  (as stated in the Lemma
\ref{LemmaLinear}), we have,
for a linear stratum:
\beq
[\tilde\Sigma]=[T_1\sim U_1][T_2\sim U_2]\dots[T_k\sim U_k]
\eeq
Therefore, the case of linear singularities is completely solved by liftings.
\subsubsection{Degeneration of non linear strata}\label{SubSubSecNonLinearStrata}
The singularities that are not linear are much harder to deal with. In such cases, the lifted stratum is not a locally
complete intersection, so we cannot use the classical intersection theory directly (this is shown by examples of $A_4,A_5$ in
section \ref{SectionA4A5}). Further liftings also do not help (as discussed in section \ref{SecDiffBlowUp}).
Instead, we use the method of degenerations.

The idea of degeneration was explained in the introduction: we want to degenerate the given stratum to strata of
higher singularities, which are nevertheless simpler to work with. In this way we get the equation for cohomology classes:
\beq\label{degen}
[\mbox{original stratum}]\times[\mbox{degenerating condition}]=[\mbox{simpler stratum}]
\eeq
The main steps of the method are:
\li To obtain the defining equations of the given stratum.
\li{To choose "the direction" of degeneration}  (which varieties we would like to obtain as a result of degeneration?).
\li{The degenerating condition}. How to define the degenerating divisor/cycle? Is the intersection of degenerating
 cycle with the original stratum transversal? What to do if it is not? Is the degeneration invertible (i.e. does the
 equation (\ref{degen}) fix the cohomology class of original stratum uniquely?)
\li Which varieties do we obtain after degeneration? What are the multiplicities?

We consider now these steps in more details, the simplest example of degeneration ($A_2\rightarrow D_4$) is considered in
section \ref{SeExamDegen}.
\paragraph{How to obtain the defining equations of equisingular stratum.}\label{ParDefiningEquationsStratum}
As will be seen in paragraph \ref{ParResultVariet}, we do not need to obtain global, covariant equations, defining
the stratum. Rather, it will always be sufficient to have equations that define the equisingular stratum of
curves with singularity at a {\it given point} (we always take: $x=(1,0,0)\in\mP_x^2$), with a {\it given
tangent line} (we always take: $l=(0,0,1)\in(\mP_l^2)^*$). Such a system of equations is obtained (in each
particular case) by elimination.

Start from a general polynomial $f=x_2+\sum_{i+j>1}a_{i,j}x_1^ix_2^j$ of degree $d$
(in local coordinates around $x=(1,0,0)$) that has $A_2$
singularity in the origin, with $x_1$-axis as a tangent line. Since the germ ($f(x)=0$) possesses the given singularity,
 it
can be brought to a preliminary form by locally analytic transformation (as the singular point and the tangent line
are fixed, we can assume that such transformation has no linear part):
\beq
x_1\ra x_1+\sum_{i+j>1} A1_{ij}x^i_1x^j_2~~~~~~x_2\ra x_2+\sum_{i+j>1} A2_{ij}x^i_1x^j_2
\eeq
 Demand that the resulting curve is in a preliminary form. This gives a set of polynomial equations (in $A1_{ij}$,
 $A2_{ij}$ and $a_{i,j}$). Eliminate from the equations all the parameters of the locally analytic
transformation ($A1_{ij}$ and $A2_{ij}$). In this way we get the basis of the locally defining ideal $I$ of the stratum
(with $x=(1,0,0),~l=(0,0,1)$).
 Two examples of this procedure (for $A_4,A_5$) are considered in section \ref{SectionA4A5}.
It is important to consider the whole basis of the ideal $I$ and not just some of its generators,
since the corresponding variety is usually not a locally complete intersection.

The aforementioned procedure can be easily done by pen\&paper in the simplest cases ($A_4,D_7$). In other
cases one can use computer programs (like Singular or Mathematica), the procedure is easily performed for every type
of singularity.

\paragraph{The degenerating condition.}\label{ParDegenCondit}
Among the generators of the locally defining ideal $I$, we consider those which are linear monomials
(i.e. give equations of the form $a_{ij}=0$). We call them: the linear part of the ideal.

The goal of degeneration is to get rid of the non-linear equations. For this we impose vanishing
of additional coefficients. To define this formally, consider the linear part of the ideal $I$.
 For example the ideal always contains the equations $a_{00}=0$, $a_{01}=0$ and $a_{10}=0$.

Let $n$ be maximal such that for every $i,j$ with $i+j<n$ one has: $a_{ij}=0$ is an equation in the basis of $I$.
Let $i_0$ be maximal such that $\forall i<i_0$ the equation $a_{i,n-i}=0$ is in the basis. Then the degeneration is by
$a_{i_0,n-i_0}\ra0$.

The cohomology class of the degenerating hypersurface is obtained in \ref{SecCohomClasesOfDegenerations}:
\beq
[degeneration~a_{pq}\ra0]=(d-p-2q)X+F+(q-p)L
\eeq

Another important issue is to check that the degeneration is invertible, i.e. equation (\ref{degen}) fixes the
cohomology class of the original stratum uniquely. From the point of view of cohomology ring
$H^*(\mP_x^2\times(\mP_l^2)^*\times\mP_f^D)$ equation (\ref{degen}) reads:
\beq
p_1(X,L,F)p_2(X,L,F)=p_3(X,L,F)
\eeq
here $p_i(X,L,F)$ are homogeneous polynomials in generators of the cohomology ring. The generators are nilpotent
($X^3=0=L^3=F^{D+1}$). Thus the solution is unique, provided: deg$(p_3)\leq D$ and $p_2$ depends in essential way on $F$.
The first condition is always satisfied, while the second means that the degeneration must involve (in essential way)
the function  $f$ (or its derivatives). In particular, conditions involving the point ($x$) and tangent line ($l$) only,
are non-invertible and will not be used.
\paragraph{The resulting varieties}\label{ParResultVariet}
After choosing the way of degeneration, we substitute the degenerating condition into the defining ideal of equations,
to check, which stratum appears as a result of degeneration. Except for the simplest cases, the stratum will be reducible
and non-reduced (with some multiplicity). One should factorize the so obtained ideal $I_{degen}$ into
its prime components (preserving the multiplicities).

The explicit calculations can be made by pen\&paper in the simplest cases only (e.g. $A_4,D_7$). However
using computer programs (Singular,Mathematica) we can easily perform all the calculations in each particular case.

An important question to consider is: to which singularity types do the so obtained ideals correspond? Here not
much can be said. The singularity types may certainly be non-semi-quasi-homogeneous or Newton degenerated. Even
worse, as we start from a complicated system of non-linear equations, the resulting ideals may correspond
to some sub-strata of the equisingular strata. And it can be quite difficult to understand the geometric
meaning of such substrata.

However, for the algorithm to work we do not need any identification with some known types, all we need is
the defining ideal to start from. Therefore this difficulty does not show up in actual calculations.

One additional point to address is the singularities of strata. Every lifting
can be thought of as being a partial resolution of a stratum. For example, the simplest lifting $\tSi(x)\ra\Si$
is the embedded blowup with its centre- the substratum of binodal curves.

The lifted strata are less singular, in some cases (linear singularities) they are even smooth
(lemma \ref{LemmaLinear}). Correspondingly, the adjacent strata appear with multiplicity one
in course of degeneration.
For non-linear strata the multiplicity is higher and should be checked explicitly in each case.
The effect of lifting $\tSi(x,l)$ is taken into account by considering curves with fixed $x,l$.

In this way a particular vertex of the degenerating tree is constructed, with outgoing edges,
decorated by multiplicities. Having completed this step, one moves to the next vertices.
\paragraph{The degenerating tree}
Now we finalize the construction of the degenerating tree.

 A particular vertex of the degenerating tree was constructed in \ref{ParResultVariet}.
We continue in this way until we get a collection of ideals generated by linear monomial
equations (i.e. of the form $a_{ij}=0$). This happens in a finite number of steps
(not bigger than the number of variables in the initial locally defining ideal).

For each so obtained ideal, we check whether it is an ideal of a linear singularity type. For this,
mark on the plane of Newton diagram the points, corresponding to the linear monomial equations.
If the envelope of this set of points is convex, then we have a convex Newton diagram,
defining some linear type.

Otherwise, one continues to degenerate, until one gets a convex Newton diagram. Again this process is finite.
So, we have constructed the degenerating tree, which is finite, has the initial types as its root
and linear singularities as its leaves.

\subsubsection{Some exceptional cases}
\paragraph{Newton degenerate singularities.}
Newton degenerate singularities are enumerated as other non-linear singularities. In this case, of coarse,
the defining non-linear conditions encoding the special degeneracy of the function must be supplied at
the very beginning.
We demonstrate this in section \ref{SecNewDegenSing}.
\paragraph{Calculations in non-universal domain}
In general one does all the calculations, assuming that $d$ (the degree of the curve) is very large. In this way, one
evades all the dangerous small-$d$ effects (coincidences on the Newton diagram, appearance of unexpected strata of
very high singularities etc.) While the calculation is done for $d$ unknown and very large, the result is valid for
$d$ in universal domain.

When $d$ is small, compared to the codimension of the singularity, we are in the non
universal domain, where the answer to the enumerative problem is non-polynomial in $d$.
Correspondingly, the calculations should be done separately for each specific value of $d$ (the
enumeration if by degeneration).
One should take into account various small degree effects.

For example, when working with the stratum of $A_k$ singularities, the dangerous additional piece is the stratum
of reducible curves with a double line as a component. Formally, it satisfies all the equations for any $A_k$.
For curves of high degrees this stratum is of huge codimension ($2d-1$),
and can be safely neglected. For small degrees (relatively to $k$) the dimension of this "wrong" stratum is
{\bf bigger} than that of the "true" stratum (of curves with isolated singularities).
Therefore the direct calculation (using lifting only) is impossible. One must use degenerations, and be
particularly careful to all the possible appearances of such strata. We consider example of such calculation
in section \ref{SecCalcInNonUnivDom}.

\subsection{Examples}\label{SecExamples}
\subsubsection{Complete intersections (ordinary multiple point)}\label{SecOrdMultPoint}
The ordinary multiple point is the singularity with the normal form $x^{p+1}_1+x^{p+1}_2$, another definition is
\mbox{jet$_p(f(x_1,x_2))=0$}. (For small $p$, these singularities correspond to: $A_1,~D_4,~X_9(=X_{1,0})$ etc.).
The defining condition here is vanishing of $f^{(p)}$ (tensor of derivatives of order $p$, in homogeneous coordinates).
The lifted variety in this case is:
\beq\label{node}
\tilde\Sigma(x)=\overline{\{(f,x)|~f^{(p)}|_x=0\}}\subset\mP^D_f\times\mP^2_x
\eeq
This variety is defined by ${p+2}\choose{2}$ transversal conditions.
To check the transversality of conditions in (\ref{node}) we note that $PGL(3)$ acts freely and transitively on $\mP^2$.
Therefore it is sufficient to check the transversality, at some particular point. We fix $x=(1,0,0)\in\mP^2$.
Then the conditions of (\ref{node}) are just linear equations in the space $\mP^D_f$ of all polynomials of the given
degree (which is at least $p+1$), so
the transversality is equivalent to linear independence. And this is immediately verified.

So the cohomology class of $\tilde\Sigma$ is just the product (in cohomology) of the classes of defining conditions,
namely (since all the conditions have the same cohomology class):
\beq
[\tilde\Sigma]=[f_{1,\dots,1}|_x=0]^{{p+2}\choose{2}}
\eeq
According to the remark at the beginning of section \ref{SectVar}, each class is:
\beq
[f_{1,\dots,1}|_x=0]=(d-p)X+F\in H^0(\mP^D_f)\otimes H^2(\mP^2_x)\oplus H^2(\mP^D_f)\otimes H^0(\mP^2_x)
\eeq
Thus:
\beq\label{nodefull}
[\tilde\Sigma]=\Big((d-p)X+F\Big)^{{p+2}\choose{2}}
\eeq
The precise action of Gysin homomorphism:
$H^{k+2}(\mP^D_f\times\mP_x^2)\stackrel{\pi_*}{\rightarrow}H^k(\mP^D_f)$ was explained in section
\ref{SecLiftings}. To obtain the cohomology class $[\Sigma]$ we just need to extract the coefficient of $X^2$. Thus,
finally we obtain the classical result:
\beq\label{nodedeg}
{\rm{deg}}(\Sigma)={{p+2\choose{p}}\choose{2}}(d-p)^2
\eeq
The ordinary multiple point is the simplest type of singularity. The specific feature that greatly simplified
the calculation, was the global transversality of the defining equations of $\tilde\Sigma$.
(In other words the lifted version is a globally complete intersection).

When the conditions are only locally transversal, the natural idea is: restrict the consideration to an affine part
of $\mP^2$. There, choose a sufficient set of transversal conditions. Now define a new projective variety, as the
projective closure (i.e. just homogenize equations). Check the situation at "infinity" and remove residual pieces.
As the simplest example we consider nodal curves.
\subsubsection{Removing residual pieces (nodal curve defined in affine coordinates)}\label{SubSubSectionNodal}
Here we want to demonstrate the process of removing residual pieces, which can arise when the defining conditions are
locally (but not globally) transversal.
Choose the affine part: $(x_0\neq0)\subset\mP^2_x$. The conditions of nodality, written in local coordinates, are: vanishing
of the two derivatives and of the function itself. Thus consider the variety:
\beq\label{nodeaff}
\Xi=\{(f,x)|f_1(x)=0=f_2(x),~f(x)=0\}\subset\mP^D_f\times\mP^2_x
\eeq
Over the affine part of $\mP^2_x$ this variety coincides with $\tilde\Sigma_{A_1}$, however at infinity ($x_0=0$) one
can expect some additional pieces. Indeed, the Euler formula for homogeneous polynomial of degree $d$
\beq\label{Euler}
\sum_{i=0}^2 x_i\di_if(x)=df(x)
\eeq
provides, that the equations of (\ref{nodeaff}), when translated to the neighbourhood of infinity, are:
\beq
x_0f_0(x)=0=f_1(x)=f_2(x)
\eeq
That is, the variety of (\ref{nodeaff}) is a union of
$\tilde\Sigma_{A_1}$ and some piece at "infinity" (when $x_0=0$), taken with multiplicity one
(since $x_0$ appears in the first degree).
So, to calculate the (co)homological class
$[\tilde\Sigma_{A_1}]$ one should subtract from $[\Xi]$ the (co)homological class of the variety defined by:
$\{x_0=0,~f_1=0=f_2\}$.

Calculating as previously: $[\Xi]=\Big((d-1)X+F\Big)^2(dx+f)$. The (co)homological class of
the piece at infinity is:  $\Big((d-1)X+F\Big)^2x$. Their difference gives the answer of (\ref{nodedeg}).

We emphasize again, that here we have dealt with the simplest case. In general the procedure of defining projective
closure of some affine variety and then removing some unnecessary pieces at infinity can be rather tricky.

Another method is to try to lift the $\tilde\Sigma$ to an even bigger space, where it is defined by a transversal set
of conditions. We consider now the simplest such case: cuspidal curves.
\subsubsection{An example of lifting. Cusp $x^p_2+x^{p+1}_1~~~~(A_2,~E_6,~W_{12})$}\label{SecCusp}
Cuspidal singularity is characterized by its point and direction $(x,k)\in\mP(T\mP^2_x)$.
It is easier to think of
the direction as given by a line (which in the plane is defined by a 1-form: $l\in(\mP^2_l)^*$), so:
\beq
\mP(T\mP^2)\simeq\{(x,l)|~l(x)=0\}\subset\mP^2_x\times(\mP^2_l)^*
\eeq
(the condition $l(x)=0$ means that the line corresponding to $l$ passes through $x$).

To write the condition of cuspidality, note that this singularity types is linear (can be brought to
its preliminary form by $PGL(3)$ transformations).
In a preliminary form one has: $jet_p{f}=x^p_2$
i.e. a power of some linear form (which correspond to the tangent line).
This condition is covariant under $PGL(3)$ transformations, therefore we define the lifted stratum as:
\beq
\tilde\Sigma=\{(x,l,f)|~l(x)=0~~f^{(p)}\sim (l\times\dots\times l)\}\subset\mP^2_x\times(\mP^2_l)^*\times\mP^D_f
\eeq
By direct check, the condition $f^{(p-1)}=f^{(p)}(x)=0$ is satisfied automatically.

The two conditions are certainly mutually transversal (e.g. $f$ appears in the second only). Thus, the total
cohomology class is the product of the classes of the conditions.
The (co)homology class of the hypersurface $l(x)=0$ is just $X+L$ (cf. equation (\ref{DemostratCohomClasses})).
The second condition (proportionality of two tensors), will occur frequently in the paper
and is considered in Appendix.
Its cohomology class is:
\beq
[f^{(2)}\sim(l\times l)]=\sum^5_{i=0}((d-2)X+F)^i(2L)^{n-i}\in H^{10}(\mP^D_f\times\mP^2_x\times(\mP^2_l)^*)
\eeq
Therefore, the total cohomology class of the lifted stratum is
\beq
[\tilde\Sigma]=(L+X)\sum^\frac{p(p+3)}{2}_{i=0}((d-p)X+F)^i(pL)^{\frac{p(p+3)}{2}-i}
\eeq
Extracting the coefficient of $X^2L^2F^{\frac{p(p+3)}{2}-3}$ we get
\beq [\Sigma]=\frac{p(p-1)(p+4)}{8}(d-p)\left(d(p^2+3p-2)-p(p^2+3p-6)\right)f^{\frac{p(p+3)}{2}-3}\in
H^{p(p+3)-6}(\mP^D_f)
\eeq
which gives:
\li{$p=2$}: $[\Sigma_{A_2}]=12(d-1)(d-2)$
\li{$p=3$}: $[\Sigma_{E_6}]=21(d-3)(4d-9)$
\li{$p=4$}: $[\Sigma_{W_{12}}]=24(d-4)(13d-44)$

For $p=2$ this recovers the well known result from \cite{Alufi}.
This approach is extremely effective when working with curves with linear-singularities.
\subsubsection{The simplest example of degeneration: $x^p_2+x^{p+1}_1\rightarrow x^{p+1}_2+x^{p+1}_1$}\label{SeExamDegen}
\parbox{13cm}
{(For $p=2$ the degeneration is $A_2\rightarrow D_4$.)
Here we enumerate cuspidal curves in one additional way, to demonstrate the degeneration method. We use the fact,
that the stratum of curves with the multiple point point ($x^{p+1}_2+x^{p+1}_1$) is adjacent to the
stratum of cuspidal curves.

To degenerate, we
demand that when bringing the cuspidal curve to a preliminary form ($\alpha x^p_2+x^{p+1}_1$), the
coefficient $\alpha$ of $x^p_2$ vanishes.
}
\begin{picture}(0,0)(-30,-65)
\mesh{0}{-90}{2}{2}{30}{75}{70}
\multiput(0,-32)(2,-2){30}{.}   \multiput(0,-62)(2,-1){30}{.}
\put(-25,-30){p+1}  \put(-20,-60){p}   \put(50,-100){p+1}
\put(5,-60){\vector(1,1){10}}
\end{picture}
\\
Thus we get jet$_p(f)=0$, which is the ordinary multiple point. Explicitly, start from the lifted stratum:
\beq
\tilde\Sigma_{x^p_2+x^{p+1}_1}(x,l)=\Big\{(x,l,f)|~f\mbox{ has a cusp at the point $x$ with the tangent line $l$}\Big\}
\eeq
The precise form of degenerating hypersurface is $V=\{\di_2^pf=0\}\subset\mP^D_f\times\mP^2_x\times(\mP^2_l)^*$.
This hypersurface is transversal to the lifted stratum in the open neighbordhood $l_2\neq 0$.
The direct check shows, that:
\beq
\tilde\Sigma_{x^p_2+x^{p+1}_1}(x,l)\cap V=\Big\{(x,l,f)|l(x)=0,~f^{(p)}=0\Big\}\cup
\Big\{\tilde\Sigma_{x^p_2+x^{p+1}_1}\cap (l^p_2=0)\Big\}
\eeq
We emphasize that though the codimension of $\Sigma_{x^{p+1}_2+x^{p+1}_1}\subset\bar\Sigma_{x^p_2+x^{p+1}_1}$ is two
(for example for $D_4$ and $A_2$), the
codimension of the lifted stratum is 1 (since
$\tilde\Sigma_{x^{p+1}_2+x^{p+1}_1}(x,l)\ra\tilde\Sigma_{x^{p+1}_2+x^{p+1}_1}(x)$ is a $\mP^1$ fibration).

Note, that the second variety is of multiplicity $p$ (not reduced). Moving it to the left,
we have for the cohomology classes:
\beq\label{CohomRingEquat}
[\tilde\Sigma(x,l)]\Big([\di_0^pf=0]-p[l_0=0]\Big)\!\!=\!\![f^{(p)}=0,l(x)=0]
\eeq
So the cohomology class of the degenerating divisor in this case is: $([V]-p[l_2=0])$. Note, that we know the
cohomology class, without writing explicit equations for the divisor. This is the typical situation during
degenerations.

The right hand side is a stratum of curves with ordinary multiple point (complete intersection)
its (co)homology class has been already
calculated.
According to the remark in section \ref{ParDegenCondit} the class of degenerating divisor is "invertible" (i.e.
the last equation has unique solution). Thus from this formula (by division) we obtain the cohomology class of
the cuspidal stratum.
So this gives a possibility to calculate the cohomology class $[\tilde\Sigma(x,l)]$.
Since the l.h.s. was also calculated, the example provides a simple check of the degeneration idea.
\section{Enumeration of linear strata}\label{SecLinearStrata}
In this section we work with Newton non-degenerate singularities only. We also assume that the degree of the curve
is sufficiently high, compared to the codimension of the singularity.

In the first part of this section we use the lifting method to calculate the cohomology classes of linear strata.
The calculation is straightforward in each particular case.
In the second part we present the method of degenerations to calculate the same things. For linear strata
it is just another check of the already obtained answers. For nonlinear strata,
however, it is the only currently working method.
\subsection{The definition and properties of linear strata}\label{SecDefinitPropertLinearStrata}
\bed
The {\bf{linear stratum}} is the equisingular stratum of (Newton-non-degenerate) singular curves that
can be brought to a preliminary form by
projective transformations only (or linear transformations in the local coordinate system centered at the
singular point). Otherwise the stratum is called nonlinear.
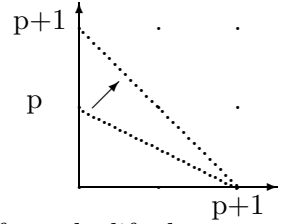
\eed
As will be shown in Lemma \ref{LemmaLinear}, the linear stratum can be lifted to a variety, defined by equations linear in function or
 its derivatives, and therefore easy to work with.
We call the singularity {\it (non)linear} iff the corresponding equisingular stratum is (non)linear.
\\
\parbox{14cm}
{The following simple property gives a useful characterization of linear strata
\bec
The equisingular stratum is linear iff every segment of the Newton diagram
of the corresponding singularity has the slope  $\frac{1}{2}\leq\mbox{tg}(\alpha)\leq2$.
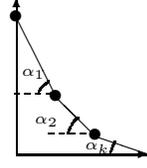
\eec
}
\begin{picture}(0,0)(-10,20)
\put(0,0){\vector(0,1){60}}\put(0,0){\vector(1,0){50}}
\put(-3,50){$\bullet$} \put(0,52){\line(1,-2){15}}
\multiput(13,23)(-4,0){4}{\line(-1,0){2}} \put(2,30){\tiny$\alpha_1$}  \qbezier(12,28)(8,25)(8,23)
\put(12,20){$\bullet$} \put(13,24){\line(1,-1){15}}
\multiput(30,8)(-4,0){5}{\line(-1,0){2}}  \put(7,12){\tiny$\alpha_2$} \qbezier(23,14)(21,13)(19,8)
\put(27,5){$\bullet$} \put(30,7){\line(3,-1){20}}
\put(26,2){\tiny$\alpha_k$} \qbezier(37,5)(36,3)(35,0)
\end{picture}\\
\bpr$\Leftarrow$
Suppose all the slopes are bounded as above and a hypersurface germ has been brought to the given
preliminary form by a chain
of locally analytic transformations. Start undoing these transformations to achieve the initial germ. Immediate
check shows that any nonlinear analytic transformation (without linear part) has no effect on the points under the Newton diagram.
Correspondingly the monomials of the initial polynomials that lie under the Newton diagram are restored by linear
transformation only. But it means that the germ could be brought to the preliminary form by linear transformations only.
\\$\Rightarrow$
Suppose at least one of the angles (in a preliminary) does not satisfy the condition $\frac{1}{2}\leq\mbox{tg}(\alpha)\leq2$. Then there exists
a quadratic shift of coordinates that changes the Newton diagram. Of course, such shift cannot be undone by linear transformations.
\epr

Another characterization of linear singularity is via the resolution tree with degrees of intersections of strict
transforms with exceptional divisors assigned to each infinitely near point (section \ref{SecTypesSingul}).
\bec
If the singularity is linear, then  there are no 3 equal consequent numbers ascribed to the vertices of the tree
\eec

The simplest class of examples of linear singularities is defined by the series: $f=x^p+y^q,~~p\leq q\leq2p$.
In general, for a given series only for a few types of singularities the strata can be linear.
In the low modality cases the singularities brought to the preliminary form by projective (linear)
transformation are (all the notations are from \cite{AGV}, see also section \ref{SecTypesSingul}):
\li{Simple singularities (no moduli)}: $A_{1\le k\le3},~~D_{4\le k\le6},~~E_{6\le k\le8}$
\li{Unimodal singularities}: $X_9(=X_{1,0}),~~J_{10}(=J_{2,0}),~~Z_{11\le k\le13},~~W_{12\le k\le 13}$
\li{Bimodal}: $Z_{1,0},~~W_{1,0},~~W_{1,1},~~W_{17},~~W_{18}$

We emphasize again, that singularities are considered up to topological transformations, this approach removes moduli.
Henceforth the moduli are ignored. For example by the singularity $X_9$ we mean singularity
defined by vanishing of all the derivatives up to (including) order 3, or (to save the words) by the line
through (4,0) and (0,4) on Newton diagram.

Most equisingular strata are nonlinear. For example if a curve has an $A_4$ point, the best we can do by projective
transformations is to bring it to the normal form of $A_3$:
\beq
a_{0,2}x^2_2+a_{2,1}x_2x^2_1+a_{4,0}x^4_1
\eeq
This quasi-homogeneous form is degenerated ($a^2_{2,1}=4a_{0,2}a_{4,0}$) and by quadratic (nonlinear!) change of coordinates
the normal form of $A_4$ is achieved.

However, inside every (equisingular) non-linear stratum there is a linear substratum: collection of points corresponding to
the curves with the same singularity that can be brought to the normal form by the {\it projective} transformations only. E.g. in
the (closure of the) stratum of $A_k$ singularities there is a substratum of singular curves that can be brought to the normal
form of $A_k$ by {\it{projective}} transformation only. We will say that these curves possess $A^{(l)}_k$ singularity, they span
the substratum $\Sigma_{A^{(l)}_k}$ of $\Sigma_{A_k}$. (If the singularity $X$ is by itself linear then of course:
$X^{(l)}\equiv X$).

Even more generally, we will consider the strata as
 \mbox{$\Sigma_{A_k}\cap\Sigma_{A^{(l)}_j}$}, $k>j>3$ ($A_k$ singularities, that can be brought to the normal form of $A_j$ by
projective/linear transformations only).
These strata will be helpful, eventhough they have no direct geometric meaning and not much importance in the general
singularity theory.
\\
\parbox{12cm}
{To distinguish between an equisingular stratum and its linear substratum, the Newton diagram of each singularity will contain two lines:
solid (corresponding to the normal form) and dotted
(corresponding to the "optimal" form to which the singularity can be brought by projective/linear transformations).}
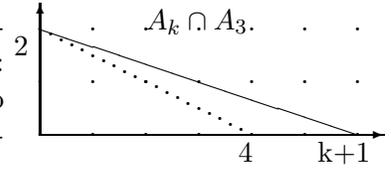
\begin{picture}(0,0)(-10,-40)                  \put(40,-10){$A_k\cap A_3$}
\mesh{0}{-50}{6}{2}{20}{130}{50}
\multiput(-1.5,-10.5)(4,-2){20}{.}
\put(0,-10){\line(3,-1){120}}   \put(-10,-20){2}   \put(75,-60){4}  \put(105,-60){k+1}
\end{picture}
\subsection{The direct approach: liftings}
The enumeration of linear strata is completed by the following
\bel\label{LemmaLinear}
Every linear stratum $\Sigma$ can be lifted to a variety $\tilde\Sigma$, which is a locally complete intersection. $\tilde\Sigma$
is a variety of specific type, defined by transversal intersection of conditions, each of them being proportionality of symmetric
multi-forms: $\Omega^{(p_i)}_i\sim \omega^{(p_i)}_i$. Correspondingly, the cohomology class is calculated as a product:
\beq
[\tilde\Sigma]=\prod_i[\Omega^{(p_i)}_i\sim \omega^{(p_i)}_i]
\eeq
\eel
The proof of the lemma is essentially the calculation of the classes by the following algorithm:
\begin{enumerate}
\item{Record the conditions of the preliminary form from the Newton diagram}
\item{Write them in a covariant way, using the one-form $l$ corresponding to the tangent line and
"resolve" the conditions using auxiliary forms, i.e. lift the stratum to a bigger space}
\item{Use the formula for the (co)homology class of condition of proportionality of two tensors
to calculate the total (co)homology class}.
\end{enumerate}
\pf

{\bf Step 1.}
To obtain the defining equations of the lifted stratum, write the conditions for the preliminary form of the singularity.
\\
\parbox{11cm}
{These are read off the Newton diagram: derivatives of $f$, corresponding to the points below the defining
(piecewise-linear)
curve, must vanish. So, one gets a set of conditions (in local coordinates) of the form:
\beq\ber
\di^i_{x_1}f=0,~~i=0,1,\dots,i_0\\\di^i_{x_1}\di_{x_2}f=0,~~i=0,1,\dots,i_1\\\dots\\\di^i_x\di^r_yf=0,~~i=0,1,\dots,i_r
\eer\eeq
}
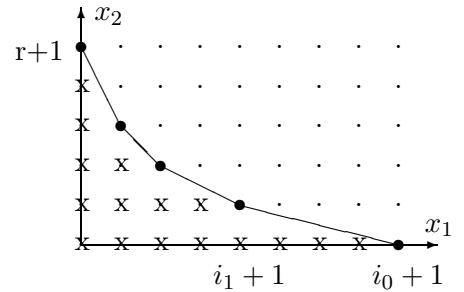
\begin{picture}(0,0)(-30,-30)
\mesh{0}{-70}{8}{5}{15}{135}{90}
\put(-25,0){r+1}           \put(50,-85){$i_1+1$}     \put(110,-85){$i_0+1$}
\put(0,5){\line(1,-2){15}}  \put(15,-25){\line(1,-1){15}}   \put(30,-40){\line(2,-1){30}} \put(60,-55){\line(4,-1){60}}
\put(0,5){\circle*{4}}    \put(15,-25){\circle*{4}}  \put(30,-40){\circle*{4}}  \put(60,-55){\circle*{4}}  \put(120,-70){\circle*{4}}
\put(5,15){$x_2$}  \put(130,-65){$x_1$}
\put(-2.5,-12){x}  \put(-2.5,-27){x} \put(-2.5,-42){x} \put(-2.5,-57){x} \put(-2.5,-72){x}
   \put(12.5,-42){x} \put(12.5,-57){x} \put(12.5,-72){x}
    \put(27.5,-57){x} \put(27.5,-72){x}
     \put(42.5,-57){x} \put(42.5,-72){x}  \put(57.5,-72){x} \put(72.5,-72){x} \put(87.5,-72){x} \put(102.5,-72){x}
\end{picture}
\\
\\
The key idea here is that the above conditions can be written as the decomposition of jets. Namely:
\beq
jet_{i_0}(f)=x_1(...),~~~jet_{i_1}(f)=x^2_1(...),~~~...jet_{i_r}(f)=x^r_1(...),~~~
\eeq
{\bf Step 2.} These last conditions are covariant under $PGL(3)$ rotations. Therefore they can be resolved,
by introducing new variables (symmetric forms) $\{B^{(i_j-1)}_j\}^k_{j=1}$:
\beq\ber
f^{(p_1)}\sim\mbox{SYM}(\underbrace{l,\dots,l}_{p_1-i_1+1},B^{(i_1-1)}_1)~~\dots~~\dots~~\dots~~
f^{(p_k)}\sim\mbox{SYM}(\underbrace{l,\dots,l}_{p_k-i_k+1},B^{(i_k-1)}_k)\\f^{(p_{k+1})}=0,~~l(x)=0
\eer\eeq
Here on the r.h.s. the one-form $l$ and the symmetric $(i_j-1)$-forms $B^{(i_j-1)}_j$ are symmetrized, the conditions mean
proportionality of tensors. Using the Euler's formula ($\sum x_i\di_i f=df$) and its consequences, the conditions can be rewritten
in the following form:
\beq\label{LiftedStratumConditions}\ber
f^{(p_1)}\sim\mbox{SYM}(\underbrace{l,\dots,l}_{p_1-i_1+1},B^{(i_1-1)}_1)\\
B^{(i_1-1)}_1(\underbrace{x,\dots,x}_{p_1-p_2})\sim\mbox{SYM}(\!\!\!\!\underbrace{l,\dots,l}_{p_2-p_1+i_1-i_2}\!\!\!\!,B^{(i_2-1)}_2)\\
\dots\\
B^{(i_k-1)}_k(\underbrace{x,\dots,x}_{p_k-p_{k+1}})=0~~~l(x)=0
\eer\eeq
The so obtained conditions are mutually transversal (e.g. $f$ appears in the first row only, $B^{(i_1-1)}_1$ in the
first and the second, etc.).
Therefore the cohomology class of the lifted variety is just the product of the classes of conditions in
(\ref{LiftedStratumConditions}). This proves the lemma.
\proofend

{\bf Step 3.}
The lemma reduces the calculation of the cohomology classes of the lifted strata to the calculation of the classes of conditions of the
type: two tensors are proportional. This problem is considered in Appendix (\ref{TensorProportionality}). The answer is as follows: if the cohomology class of an
element\footnote{By a cohomology class of an element of a tensor we mean the cohomology class of a hypersurface defined by vanishing
of this element. We assume that at least one of the tensors is generic, in particular non-degenerate}
of a k-indexed tensor $a^{(k)}$ is $A$, while that of $b^{(k)}$ is $B$, then:
\beq
[a^{(k)}\sim b^{(k)}]=\sum^{{k+2\choose{2}}-1}_{i=0}A^iB^{{k+2\choose{2}}-1-i}
\eeq
Many examples of the enumeration are given in Appendix. We solve there cases of singularities with normal forms:
$x^p_2+x^{p+k}_1$~$0\le k\le3$~~$(A_{k\le3},E_6,E_8,J_{10},W_{12},W_{1,0},W_{18}),$
$x^py+y^q$,~$\frac{p}{2}<q-1\leq2p$~~$(D_{k\le6},E_7,Z_{11},Z_{13},Z_{1,0},W_{13},W_{17})$ and some others.
The method works equally well for the linear non quasi-homogeneous singularities (defined by two or several lines
on the Newton diagram). We solve in Appendix the cases of: $X_{1,p},~0<p\leq2$ and $W_{1,p},~0<p\leq2$.
\subsection{Degenerations}
The process of degeneration was described in section \ref{SubSubSecNonLinearStrata}. In general the related
calculations are rather involved. However, for linear singularities, they are particularly simple. One sees
directly from the Newton diagram the degenerating direction and the final stratum.
And after each degeneration the stratum is irreducible (and reduced).

Since the simplest type of singularity is ordinary multiple point, we would like to degenerate the (linear) equisingular
stratum to the stratum of some ordinary multiple point.

So, the goal of the degeneration process (in the case of linear singularities) is to deform the piecewise linear
curve, defining the singularity on the Newton diagram, to a line, and to rotate the line clockwise.
At each step of degeneration one imposes additional condition: the coefficient of some monomial in a preliminary
form vanishes. The precise degeneration is discussed in \ref{ParDegenCondit}.
\section{Enumeration of nonlinear strata}
In this section we treat the non-linear strata.
Except for the simplest cases at the beginning of each series (e.g. $A_1,A_2,A_3$ for $A_k$), the normal form of the
singularity cannot be achieved by projective transformation. So, the (lifted) stratum of curves (with singularity
at a prescribed point, with a prescribed tangent line) cannot be defined
(even locally) by equations linear in function or its derivatives. The (non-linear) equations are complicated and difficult to
deal with.

Due to this reason, the direct approach does not work (we consider examples of $A_4,A_5$ in the first subsection).
The further lifting faces severe difficulties and gives no results (this is demonstrated for the
case of $A_k$).
The only currently working method is "degeneration", it was described in section \ref{SubSubSecNonLinearStrata}.
It provides a way to calculate the cohomology classes of non-linear strata.

Usually there are many different ways to degenerate the given stratum. The most efficient are the ones that push
the given stratum to a linear stratum in the smallest number of steps and involve the least number of other strata.

The degeneration methods work equally well for any type of singularity.
First we consider the Newton-non-degenerate singularities. In subsection \ref{SecNewDegenSing}
we consider the Newton-degenerate singularities and solve some cases.
\subsection{Impossibility of the direct approach}\label{SecImpossibilityOfTheDirectApproach}
\subsubsection{ Examples of $A_4,~A_5$}\label{SectionA4A5}
Suppose that by projective transformation
the curve (with $A_4$ singularity) is brought to the form (in the local coordinates):
\beq\label{A4local}
a_{0,2}x^2_2+a_{2,1}x_2x^2_1+a_{4,0}x^4_1+\mbox{higher~terms}
\eeq
By one quadratic shift ($x_2\rightarrow x_2+\beta x^2_1$) we should be able "to kill" the monomials $x_1x^2_2$ and $x^4_1$.
This
means that the quasi-homogeneous form in (\ref{A4local}) is degenerate, in particular the (local) condition
for $A_4$ is:
\beq\label{A4coeff}
4a_{0,2}a_{4,0}=a_{2,1}^2
\eeq
Similar considerations give the additional condition for $A_5$:
\beq\label{A5coeff}
a_{1,2}a_{4,0}+a_{5,0}a_{0,2}=\frac{a_{2,1}a_{3,1}}{2}
\eeq
The first task is to make these conditions covariant under $PGL(3)$. For this lift $\Sigma$ to:
$\tilde\Sigma\subset\mP^2_x\times\mP^2_v\times(\mP^2_l)^*\times\mP^D_f$. Here $x$ is the singular point,
$v$-another point in the plane, the tangent line $l$ passes through the two points. The equations
(\ref{A4coeff},\ref{A5coeff}) are written for the case:
\beq
x=(1,0,0)~~v=(0,0,1)~~l=(0,1,0)
\eeq
In terms of derivatives (\ref{A4coeff}) reads: $\frac{1}{3}f^{(2)}_{22}f^{(4)}_{1111}=(f^{(3)}_{112})^2$. The prescription to convert this into
covariant expression is as follows:
\begin{enumerate}
\item Every time index "1" appears; contract the derivative with v
\item Every index "2" should be replaced by a free index
\end{enumerate}
In this way, e.g. the equations for $A_4,~A_5$ cases become:
\beq\ber
\frac{1}{3}f^{(2)}_{ij}f^{(4)}(v,v,v,v)=f^{(3)}(v,v)_if^{(3)}(v,v)_j\\
f^{(3)}(v)_{ij}f^{(4)}(v,v,v,v)+f^{(2)}_{ij}f^{(5)}(v,v,v,v,v)=\frac{f^{(3)}(v,v)_if^{(4)}(v,v,v)_j}{2}
\eer\eeq
Even though the covariant conditions have been obtained, we face severe difficulties with actual cohomology computation: the
total set of conditions is certainly not transversal (even locally).
%
%
%
%
%
%
\subsubsection{The difficulties of further liftings}\label{SecDiffBlowUp}
The natural way to attack the $A_k$ stratum is to lift it further. In this approach, however,
we meet the following obstacle.

By definition, a curve has $A_k$ singularity at the origin if the $k$-jet of the function can be brought (by polynomial
transformations) to the form: $jet_k(f(x_1,x_2))=x^2_2$.

Now, go back to the original form of the function, i.e. undo all the transformations. Then we arrive at the condition:
\beq
jet_k(f(x,y))=jet_k(g^2(x,y))
\eeq
Here $g$ is a polynomial (vanishing at the singular point) of degree not more that $k-1$. (Terms of higher order do
not enter the $k$-jet and therefore can be omitted). This condition can be written as proportionality of
tensors of $k$'th derivatives:
$f^{(k)}\sim(g^2)^{(k)}$. Thus the lifted stratum of curves with $A_k$ point is ({\bf naively}) defined as:
\beq\label{wrongAk}
\overline{\{(x,f,g)|~f^{(k)}\sim(g^2)^{(k)},~g(x)=0,~~g~\mbox{is smooth}\}}\subset\mP^2_x\times\mP^{D_f}_f\times\mP^{D_g}_g
\eeq
This is, however, not the true stratum. For example the closure of $\Sigma_{A_k}$ contains the $\Sigma_{D_{k+1}}$ stratum (in
codimension 1). The stratum defined in (\ref{wrongAk}) does not contain $\tilde\Sigma_{D_{k+1}}$. Indeed, if the curve has a
$D_k$ point then the tensor of second derivatives at the point vanishes: $f^{(2)}=0$. This causes the gradient of $g$ to vanish:
$\nabla g=0$. But then $(g^2)$ has a quadruple point with pairwise coinciding tangent lines. So the closure of the
stratum, defined in
(\ref{wrongAk}), does not contain curves with $D_{\tilde{k}}$ (for any $\tilde{k}$), neither $E_{\tilde{k}}$, etc. This stratum
can be thought of as a wrong compactification of the pure $A_k$ stratum, with many substrata shrank.

Another problem is the actual calculation of the cohomology class: even the class of the above "wrong" stratum cannot be
calculated.
The problem is in the open condition of smoothness of $g$. It cannot be taken into account by subtracting a cohomology class
of unnecessary piece of singular $g$, because this piece has dimension higher than the stratum (\ref{wrongAk}). In fact, if $g$ has
 a multiple point (so that $g^{[\frac{k}{2}+1]}=0$), then the condition $f^{(k)}\sim(g^2)^{(k)}$ is automatically satisfied, so we
 have a stratum of codimension: $[\frac{k}{2}]+3\choose{2}$, while the codimension of variety in (\ref{wrongAk}) is $k+2\choose{2}$.
\\
\\
Therefore the cohomology classes of non-linear strata cannot be calculated by direct approach, neither by a suitable lifting.
Instead, we will use a particular kind of degeneration (linearization) to convert the nonlinear equations to linear ones.
As a preparation we generalize the situation of $A_4,~A_5$, i.e. we discuss the "first" nonlinear equation in each series.
\subsection{The process of degeneration.}\label{SecTheDegen}
In the nonlinear case the goal of degeneration is to get rid of nonlinear equations. The
general algorithm is described in \ref{SubSubSecNonLinearStrata}. Here we mainly illustrate it for the case of $A_k$.
There are (at least) two systematic ways to implement the degeneration:
\\
\parbox{12cm}
{\li{\bf{The clockwise rotation}.} Push the line to increase its slope\footnotemark[\value{footnote}],
 so that one arrive at a diagram with a steeper line
(and this case is easier to enumerate). As an example consider $A_4$. In this case to push the line clockwise means to
demand that in the normal form $a_{0,2}x^2_2+a_{2,1}x_2x^2_1+a_{4,0}x^4_1$ the coefficient $a_{0,2}$ vanishes,
so naively the transitions is:
$A_4\rightarrow D_5$.}
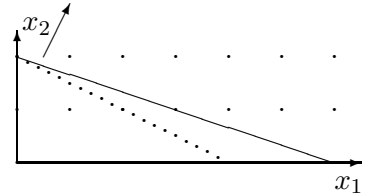
\begin{picture}(0,0)(-10,-30)        \put(10,-10){\vector(1,2){10}}
\mesh{0}{-50}{6}{2}{20}{130}{50} \multiput(-1.5,-10.5)(4,-2){20}{.} \put(0,-10){\line(3,-1){120}}
\put(2,0){$x_2$} \put(120,-60){$x_1$}
\end{picture}
\footnotetext[\value{footnote}]{as in the examples of Section 3: $x^{p+1}+y^{p+2}\rightarrow x^{p+2}+y^{p+2}$}
\\
However, from the defining (nonlinear) equation for $A_4$ ~($4a_{0,2}a_{4,0}=a_{2,1}^2$) one has:
$a_{0,2}=0=a_{2,1}^2$, i.e. the stratum $D_5$ is adjacent to $A_4$ with multiplicity 2. So, the actual transition is:
\beq
A_4\rightarrow 2D_5
\eeq
Similar considerations for $A_5$ give: under $a_{0,2}\rightarrow0$ the transition is:
\beq
A_5\rightarrow 2(D_5\cap A_5)=2D_6\cup 2E_6
\eeq
\li{\bf Linearization.} We remind that  the singularity is described on the Newton diagram by two lines:
solid (corresponding to the normal form) and dotted (corresponding to the optimal form to which we can bring
the singularity by projective transformation only).

\hspace{-0.5cm}\parbox{11cm}
{The goal of the linearization is to push the dotted line to the solid line and finally to merge them
(so that the linear stratum is obtained).
In the case of $A_4$ this means, that by linear transformation we will be able to bring the curve to the normal form, in particular the
degeneration in this case is: $A_4\rightarrow A^{(l)}_4$, i.e. we should demand $a_{4,0}\rightarrow0$.}
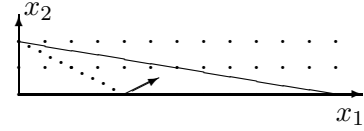
\begin{picture}(0,0)(-20,-30)        \put(40,-30){\vector(2,1){13}}
\mesh{0}{-30}{12}{2}{10}{130}{30} \multiput(-1.5,-10.5)(4,-2){10}{.} \put(0,-10){\line(6,-1){120}}
\put(2,0){$x_2$} \put(120,-40){$x_1$}
\end{picture}
\\
As in the previous case since the equation is nonlinear the degenerated stratum enters with a multiplicity.
The transition is:
\beq
A_4\rightarrow 2A^{(l)}_4
\eeq
For $A_5$ we have:
\beq
A_4\rightarrow 2(A^{(l)}_4\cap A_5)=2A^{(l)}_5\cup 2D_6
\eeq

In both cases (clockwise rotation and linearization) the key step is to recognize the strata appearing after
degeneration and their multiplicities. This is done by the analysis of
locally defining ideal of the singularity \ref{ParDefiningEquationsStratum}.

As the cases $A_4,A_5$  are considered above we start from
\bex
\li $\bf A_6$  Here the degeneration by clockwise rotation ($a_{2,0}\ra0$) results in the degeneration tree:
\beq
A_6\ra 2D_7+3E_7+2X_9,~~~ D_7\ra 2E_8
\eeq
While $E_7$ and $X_9$ are linear strata (with the cohomology classes directly calculated), the stratum $D_7$
is non-linear and should be degenerated further (again by clockwise rotation: $a_{2,1}\ra0$).

This degeneration tree gives a set of linear equations from which the class of $\tSi_{A_6}$ is restored uniquely
(the final answer in the introduction).
\li $\bf A_7$ Again, the clockwise rotation ($a_{2,0}\ra0$) gives:
\beq
A_7\ra 2D_8+3E_8+2\tilde{X}_9,~~~D_8\ra 4J_{10}+2X_{1,1}
\eeq
Here $\tilde{X}_9$ denotes the degenerated version of $X_9$,
with an additional condition: $a^2_{13}+8a^2_{04}a_{31}=4a_{04}a_{13}a_{22}$. This condition corresponds to
a specific modulus.
It is degenerated to
a linear singularity by $\tilde{X}_9\stackrel{a_{04}\ra0}{\longrightarrow}3X_{1,1}$
\li $\bf D_7$ Immediate check gives: $D_7\stackrel{a_{21}\ra0}{\longrightarrow}2E_8$
\li $\bf D_8$ Similarly $D_8\stackrel{a_{21}\ra0}{\longrightarrow} 4J_{10}+2X_{1,1}$ with both resulting types being linear.
\eex
\subsection{Newton-degenerate singularities}\label{SecNewDegenSing}\label{SecAppNewtDegen}
Up to now we have considered the Newton-non-degenerate singularities, their normal form is fixed by the Newton diagram. For
Newton-degenerate singularities, the monomials of the normal form satisfy some additional conditions, except for those read
off the diagram. Namely, some quasi-homogeneous forms corresponding to the straight segments of the piece-wise linear
Newton diagram are degenerate.

The simplest examples of Newton-degenerate singularities are:
$W_{1,2q-1\ge1}^\sharp,$$(x^2_2+x^3_1)^2+x_2x^{4+q}_1$ and $W_{1,2q\ge2}^\sharp,$$(x^2_2+x^3_1)^2+x^2_2x^{3+q}_1$
with Milnor number $\mu=15+q$.

Instead of developing the general theory, we consider a subclass of such singularities, to
demonstrate applicability of the method.
\\
Consider the singularities with (quasi-homogeneous) normal form:\\
\parbox{12cm}{
\beq\label{EqNewtDegenSing}\ber
f(x_1,x_2)=\sum^k_{i=0} a_{k-i,i}x^{(k-i)q}_1x^{ip}_2\\
{\rm{gcd}}(p,q)=1,~~~p<q<2p,~~k>1
\eer\eeq
Since $k>1$ the form contains at least three monomials. For generic
coefficients $\{a_{i,j}\}$ this is a linear singularity, the case completely solved by liftings in subsection 3.1.
}
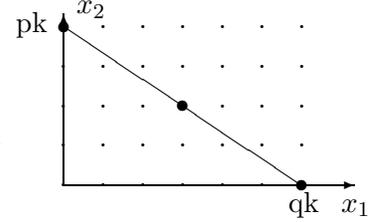
\begin{picture}(0,0)(-20,0)
\mesh{0}{-20}{6}{4}{15}{110}{65}
\put(0,40){\line(3,-2){90}}
\put(-18,38){pk} \put(5,45){$x_2$}   \put(85,-30){qk} \put(105,-30){$x_1$} \put(0,40){\circle*{4}}  \put(90,-20){\circle*{4}}
\put(45,10){\circle*{4}}
\end{picture}

For degenerate quasi-homogeneous forms the singularity will be Newton-degenerate with
the defining equations
non-linear in the derivatives of $f(x,y)$.

 As with any non-linear stratum, the strategy is to degenerate the stratum to higher linear strata. In this case the most efficient
type of degeneration is the "clockwise rotation", described at the beginning of subsection 4.2. Namely,
demand that in a preliminary
form the coefficient of $a_{0,pk}y^{kp}$ vanishes. Since the form is degenerate, this causes other monomials to vanish, and we arrive
at a "less-nonlinear" stratum (i.e. the number of nonlinear equations decreases). Degenerate further, until one gets a linear stratum.

We consider here the simplest example when the quasi-homogeneous form consists of three monomials (i.e. $k=2$).
In this case the nonlinear equation
is:
\beq
4a_{2p,0}a_{0,2q}=a_{q,p}^2
\eeq
The degeneration $(a_{2p,0}=0)$ gives (twice) the L-singularity which (in general) is not quasi-homogeneous (the Newton
diagram contains several segments). The classes of such strata are calculated immediately, however we cannot write one
general formula for arbitrary $p,q$ (it is absolutely cumbersome). Rather, it is better to solve case by case, we consider just two
simplest cases. (Denote the initial and the final strata by $\tilde\Sigma_0,\tilde\Sigma_1$.)
\begin{itemize}
\item{$\mathbf q=p+1:$} Here the degenerated singularity has preliminary form: $x^{2p}_2x_1+x^p_2x^{p+1}_1+x^{2(p+1)}_1$.
The formula for (co)homology classes is:
\beq
[\tilde\Sigma_0]=
2\frac{[\tilde\Sigma_1]}{[f_{\underbrace{(1,\dots,1)}_{2p}}=0]-p[l_1=0]}
\eeq
\item{$\mathbf q=p+2:$} Here the degenerated singularity has preliminary form:
$x^{2p}_2x_1+x^{2p-[\frac{p}{2}]}_1x^{[\frac{p}{2}]+1}+x^{2(p+2)}_1$.
The formula for (co)homology classes is:
\beq
[\tilde\Sigma_0]=
2\frac{[\tilde\Sigma_1]}{[f_{\underbrace{(1,\dots,1)}_{2p}}=0]-p[l_1=0]}
\eeq
\end{itemize}
Here $\tilde\Sigma$ is the (lifted) variety of curves with prescribed singularity at the given point, with the given
tangent line. The cohomology classes of the varieties on the r.h.s. of the formulae are calculated in Appendix.
By quotient we mean the polynomial obtained after division (by construction the numerator is divisible by denominator).

In general, the diagram of the Newton-degenerate singularity consists of segments corresponding to (degenerate) quasi-homogeneous forms.
In the same way as explained above, we degenerate each form, until we get a collection of Newton-non-degenerate linear singularities.
\section{Singular curves in the region of non-universality}\label{SecNonUniversalDomain}
As was noted in the introduction, the linearization algorithm applies for sufficiently high degrees of curve (compared to the codimension
of singularity). The final answer, however, is valid for a broader range of degrees (for a given codimension of singularity), from the
universality \cite{Kaz1},\cite{Liu} it follows that the answer is valid for the range of universality.
We obtain here two sufficient conditions for the validity of results, and explain how to calculate in the non-universal
domain.
\subsection{The domain of universality}
The method certainly works as far as we consider the stratum whose generic point corresponds to an
irreducible curve. The substrata corresponding to reducible curves may accidentally appear in the calculations, and cause
wrong result.

If the curve (with some singularity) decomposes, its components have to intersect with some tangency.
The direct check of the possible
cases of reducibility shows, that among the strata of reducible curves, the stratum with the biggest dimension
(and therefore the "most dangerous") is the stratum of
reducible curves, with a double line component. Its codimension is:
\beq
2d-1
\eeq
So the algorithm can give wrong results (and we're out of universality domain)
if the codimension of the singularity is higher than $(2d-2)$. Even if this does not happen, one could expect that in
some low degree cases various specific coincidences occur, and the bound ($\mbox{codim}<2d-1$) is not sufficient.
Later we bring some arguments that it is sufficient.
First, we derive an auxiliary criterion.

The idea of the linearization algorithm gives us a simple sufficient condition.
\bel
The linearization algorithm works and gives the correct results, if the
order of determinacy is not higher than ($d+1$).
\eel
So, for example, sufficient conditions in the simplest cases are:
\beq
A_k:k\le d,~~~D_k:k\le d+1,~~~E_{6k}:3k\le d,~~~E_{6k+1}:3k\le d,~~~E_{6k+2}:3k+1\le d,
\bet a point of\\multiplicity $p$\eet\!\!\!\!:p\le d
\eeq
\pf
First note that the statement is immediate for linear singularities. Otherwise, consider degenerations.

The algorithm degenerates the stratum to the strata of higher codimensions. At the end we represent some
combination of linear strata as a degeneration of the original stratum (or vice versa we represent the original stratum as
a deformation of the linear strata.) So the algorithm works if the singularities at the end lie in the universality
domain. The order of determinacy of the strata at the end is not bigger than that of the original stratum.
Since the strata at the end are linear (and for linear strata the lemma is immediate), we get the final answer.
\proofend

The conditions of the above lemma provide that the algorithm works, produces correct results at each step, and the final
result is also true. This bound is effective for linear singularities, but not for general ones.
We can try and improve this bound (or rather derive another one) by using degenerations.
\\{\bf Conjecture}\\
The results of the linearization algorithm are valid for the degree $d$ of the curve, and the codimension of
the singularity, satisfying $2d-1>$\mbox{codim}
\\{\bf A sketch of proof.} (We believe that the following arguments can be made rigorous.)\\
We can assume that the substratum of double line does not interfere (by the condition of the lemma).
Further, it is immediate that for
a given codimension, the singularity of type $A_k$ is the "most dangerous", since for a given codimension
it is the type with the highest order of determinacy. So it is sufficient to check the case of $A_k$.

Use a special kind of degeneration. Note that the variables ($a_{ij}$) appearing in the locally
defining ideal are naturally bi-graded: $w(a_{ij})=(i,j)$, and each equation in the ideal has terms of a constant
bi-grading. In addition each equation has constant degree in $a_{ij}$.
 Therefore, if the variable $a_{0k}$ is the one with the highest second weight, in
 any equation it is multiplied by "a partner": $a_{ij}$.
For example, for $A_k$ the higher monomial is $a_{0,k}$, while "its partner" is $a_{2,0}$.
It is this partner that is sent to zero at each degenerating step. For semi-quasi-homogeneous
singularities this coincides with the usual clockwise rotation.

So, we build a special degeneration tree and from this tree we extract only the types with the highest
order of determinacy. This gives the "most dangerous branch", to which we restrict the consideration.
The initial part of this branch is (for $k$ high enough):
\beq
A_k\ra D_{k+1}\ra J_{2,k-8}\ra \Big\{J_{i,k-ri}\ra\dots\Big\}_{i<[\frac{k}{4}]}\ra X_{3,*}\ra\dots
\eeq
The branch continues until a linear singularity is reached.
Two important facts are:
\li For $k>4$ each degeneration reduces the order of determinacy by at least one.
\li For $k>15$ after $M$ degenerating steps the multiplicity of the type is not bigger than $2+\frac{M-2}{2}$

Suppose now that after $M$ steps we have reached a linear type. Let $m$ and $o.d.$ be the
multiplicity and order of determinacy of this type. As it is linear we have: $2m\ge o.d.-1$. On the other side,
from the estimation above we get: $m<2+\frac{M-2}{2}$, while $o.d.\le k+1-M$.

Thus we get $o.d.\le Min(k+1-M,M+3)$.
The parameters $(o.d.,d)$ are in the universal domain provided: $o.d.-1\le d$. Thus, combining
all the inequalities above we get a sufficient condition of the universal domain for $k>15$
\beq
k<2d-1
\eeq
For the cases $2\le k\le 15$ one verifies the sufficiency of this bound by constructing
 the "most dangerous branch" explicitly. Some of the cases are:
\beq\ber
A_4\ra D_5,~~A_5\ra D_6,~~A_6\ra D_7\rightarrow E_8:~(o.d.=5)~~,A_7\ra D_8\rightarrow J_{10}:~(o.d.=6)\\
A_8\rightarrow D_{9}\rightarrow J_{2,1}\rightarrow X_{2,2}:~(o.d.=6),~~
A_9\rightarrow D_{10}\rightarrow J_{2,2}\rightarrow X_{2,3}:~(o.d.=7)\\
A_{10}\rightarrow D_{11}\rightarrow J_{2,3}\rightarrow E_{14}\rightarrow (x^3_2x_1+\dots+x^7_1):~(o.d.=7)
\eer\eeq
So the condition $2d-1>k$ is sufficient for the types $A_k$. Since these are "the most dangerous" types
this proves the lemma.
\proofend

The two lemmas give two sufficient conditions on the range of validity of final results.
In some simplest cases they give:
\beq\ber
A_k:~~k<2d-1~~~~D_k:~~k<2d-1~~~~E_{6k}:~~k<\mbox{max}(\frac{d}{3},\frac{2d-2}{5})~~~
E_{6k+1}:~~k<\mbox{max}(\frac{d}{3},\frac{2d-3}{5})\\
E_{6k+2}:~~k<\mbox{max}(\frac{d-1}{3},\frac{2d-4}{5})
\eer\eeq
For linear strata, the first bound (on the total degree of the boundary monomial) is a better estimation than the second,
e.g. for an ordinary p-tuple point the first bound gives $p<d$, while the second ${p\choose{2}}<2d-1$. The second
bound is better for (very) non-linear strata, e.g. for an $A_k$ point.
\subsection{Calculations in the non-universal domain}\label{SecCalcInNonUnivDom}
As always the technique of degenerations works. Namely, consider:
\beq
\tilde\Sigma=\overline{\{(x,l,f)|~f~\mbox{has~isolated}~A_k~\mbox{at}~x~\mbox{with}~l~\mbox{as a tangent line}\}}
\eeq
Degenerate this stratum by "backward rotation". After several steps one get either linear strata, or reducible curves. In
both cases the answer is obtained immediately.

As the simplest example we consider quartics with $A_7$ singularity. In this case there are (at least) two ways to obtain the answer:
by degenerations and by the classical enumeration.

The genus inequality restricts such a quartic to be reducible, so we have two (smooth) conics maximally
tangent (degree of tangency=3, degree of intersection=4). The additional (bad) stratum is now:
curves that consist of double line and (arbitrary) conic. These curves (formally) satisfy the conditions of $A_k$ for any $k$.
The dimension of both strata is: 7.
\subsubsection{Calculation via degenerations.}
Consider the defining equations for the (isolated) $A_7$ singularity of a quartic. In local coordinates
(singularity at the origin with tangent line $x_1$-axis)
the equations generating the ideal are:
\beq\ber
4a_{0,2}a_{4,0}=a_{2,1}^2,~~~2a_{1,2}a_{4,0}=a_{2,1}a_{3,1},~~~2a_{3,1}a_{0,2}=a_{2,1}a_{1,2},
~~2a_{1,3}a_{0,2}=a_{0,3}a_{1,2},~~a_{2,1}a_{1,3}=a_{3,1}a_{0,3}
\\
2a_{2,1}a_{2,2}=a_{3,1}a_{2,2}+4a_{0,3}a_{4,0},~~8a^2_{4,0}a_{1,3}=-a^3_{3,1}+4a_{3,1}a_{4,0}a_{2,2},~~
a^2_{1,2}=4a_{0,2}a_{2,2}-2a_{0,3}a_{2,1}
\eer\eeq
Degenerate the stratum by $a_{0,2}=0$. The resulting stratum is the union of two non-reduced schemes:
\li $a_{0,2}=a_{2,1}=a_{0,3}=a_{1,2}=0=...$ is $X_9$ with a specific condition on the modulus
$8a_{4,0}^2a_{1,3}=a_{3,1}(-a_{3,1}^2+4a_{4,0}a_{2,2})$. This scheme enters with multiplicity 2.
\li $a_{0,2}=a_{2,1}=a_{1,2}=a_{4,0}=a_{3,1}=0$: Conic, intersecting double line at the origin.
This scheme enters with multiplicity 1.

The class of the second stratum is readily obtained. To calculate the class of the first, intersect it with
 $a_{0,4}=0$.
Then one get $a_{3,1}^3=0$ (quadruple point with two coinciding lines).
As a result of the calculation, one has:
\beq
[\tilde\Sigma_{A_7,d=4}]=(L+X)\left(7F^{10}+5(7L+6X)F^9+18(14L^2+9X^2)F^8+504LX^2F^7\right)
\eeq
In particular, the degree of the stratum is $504$.
\subsubsection{Classical calculation}
We use here the fact that the quartic must be reducible. Namely, we want to calculate the number of
maximally tangent conics passing through 7 generic points. Two cases are possible:
\begin{enumerate}
\item{\bf{7=5+2}}. One conic passes through the 5 chosen points (and is fixed by this condition), the second passes through
the 2 remaining
points (this leaves 3 dimensional linear system). The number of such conics should be multiplied by $7\choose{2}$.
\item{\bf{7=4+3}}. One conic passes through the 4 chosen points (one-dimensional linear system), another through the 3 points
(2 dimensional linear system). The combinatorial factor in this case is $7\choose{3}$
\end{enumerate}
In both cases we have a problem of enumeration of curves with prescribed tangency. The answers are obtained in Appendix,
formula (\ref{ConicMaxTangent}). So we have:
$4{7\choose{2}}+12{7\choose{3}}$.
So, the number of maximally tangent conics, passing through the 7 generic points is
\beq
\mbox{deg}=\Bigg(4{7\choose{2}}+12{7\choose{3}}\Bigg)=504
\eeq
Alternatively, one can consider the stratum of reducible quartics which are maximally tangent conics
\beq
\tSi=\{(x,l,C_1,C_2,f)|f^{(4)}\sim C_1C_2,~~(x,l,C_1,C_2)\in \tSi_{2,2}\}
\eeq
where the class of $\tSi_{2,2}$ is obtained in Appendix. An important point is that now the projection $\tSi\ra\Si$
given by $(f,C_1,C_2)\ra f$ is $2:1$, so the cohomology class of the total variety should be divided by 2.
\appendix
\section{Appendix}
\subsection{Cohomology classes of some particular cycles}
\subsubsection{Cohomology class of a diagonal subvariety}
\label{TensorProportionality}
In the paper we often meet the conditions of proportionality of two (symmetric) tensors $a\sim b$. It defines
the diagonal subvariety:
\beq
\Delta=\{(a,b)|a\sim b\}\subset\mP_a^N\times\mP_b^N
\eeq
The cohomology class of such a condition is easily expressed (e.g. it is given in \cite[Chapter 14]{Ful}).
Namely, let $A,B$ be the generators of the cohomology rings (e.g. $H^*(\mP_a^N)=\frac{\mZ[A]}{A^{N+1}}$).
Then:
\beq
[\Delta]=\sum_{i=0}^n A^iB^{N-i}
\eeq
\subsubsection{Cohomology classes of degenerations}\label{SecCohomClasesOfDegenerations}
\parbox{12.5cm}
{Here we calculate the cohomology classes of degenerating condition $a_{pq}\ra0$ for the stratum
of curves whose locally defining ideal includes the equations:
\beq
\rm{if}~\Big\{\ber i+j<p+q\\i\le p\eer~\rm{then}~a_{ij}=0,~~~~\rm{if}~i+j=p+q~\rm{and}~i<p~\rm{then}~a_{ij}=0
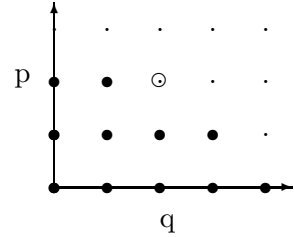
\eeq
}
\begin{picture}(0,0)(-20,-110)
\mesh{0}{-130}{4}{3}{20}{90}{70}
\put(-15,-90){p}  \put(40,-145){q}    \put(39.5,-89.5){\circle{5}}
 \put(0,-90){\circle*{4}} \put(20,-90){\circle*{4}}   \put(0,-110){\circle*{4}}  \put(0,-130){\circle*{4}}
\put(20,-110){\circle*{4}}  \put(40,-110){\circle*{4}} \put(0,-110){\circle*{4}}  \put(60,-110){\circle*{4}}
\put(20,-130){\circle*{4}} \put(40,-130){\circle*{4}}\put(0,-130){\circle*{4}}\put(60,-130){\circle*{4}}\put(80,-130){\circle*{4}}
\end{picture}
\\\\
We start from a lifted stratum $\tSi_1(x,l)$ of
curves with singularity at the point $x$, with the tangent line $l$.
 The above conditions force
(cf. section \ref{SecLinearStrata}):
\beq
f^{(p+q)}\sim\rm{SYM}(\underbrace{l..l}_{p}A^{(q)}),~~~~l(x)=0
\eeq
We lift the stratum further, by choosing
another point on the tangent line $y\in l$. Consider now the hypersurface defined by
$f^{(p+q)}(\underbrace{y..y}_{q})_{0..0}=0$. For $x\ne y$ and $l^p_0\ne0$ it corresponds precisely
to the condition $a_{pq}\ra0$, since when intersected with the above equations it gives:
\beq
f^{(p+q)}\sim\rm{SYM}(\underbrace{l..l}_{p}A^{(q)}),~~~~l(x)=0,~~A^{(q)}=\rm{SYM}(l,A^{(q-1)})
\eeq
Therefore we have:
\beq
[\tSi_1(x,l,y)][f^{(p+q)}(\underbrace{y..y}_{q})_{0..0}=0]=[\tSi_2(x,l,y)]+q[\tSi_1(x,l,y)\cap (x=y)]+[\tSi_1(x,l,y)\cap (l^p_0=0)]
\eeq
Here the residual piece over $x=y$ enters with multiplicity $q$ (by direct check). Note that
$[\tSi_i(x,l,y)]=[\tSi_i(x,l)](L+Y)$, while the cohomology class of the degenerating hypersurface
is $(d-p-q)X+F+qY$. Then we get:
\beq
[\tSi_1(x,l)](L+Y)((d-p-q)X+F+qY-pL)=[\tSi_2(x,l)](L+Y)+q(X^2+XY+Y^2)[\tSi_1(x,l)]
\eeq
(here we use the cohomology class of the diagonal). Note that in this expression the term $Y^2$ cancel, as it should be,
since the map $\tSi_1(x,l,y)\ra\tSi_1(x,l)$ is a $\mP^1$ fibration.

Finally, extracting the coefficient of $Y$ we get:
\beq
[\tSi_1(x,l)]\Big((d-p-2q)X+F+(q-p)L\Big)=[\tSi_2(x,l)]
\eeq
which gives the cohomology class of the degenerating divisor.

Note, that we have written the needed cohomology class, without writing explicitly the complete
defining ideal of the divisor.
\subsection{Some explicit calculations}
\subsubsection{For linear strata}\label{SecExplicitCalculationLinearStrata}
Here we perform explicit calculations for various linear strata. The method is described in section 3.

\paragraph{The series $ x^p_2+x^{p+1}_1~~~~(A_2,~E_6,~W_{12})$}
This singularity series is defined by vanishing of all the derivatives up to $p$, and the tensor $f^{(p)}$
satisfies a condition similar to the cusp case: $f^{(p)}\wedge l=0$. Hence we define:
\beq
\tilde\Sigma=\{(x,l,f)|~l(x)=0~~f^{(p)}\sim (l\times\dots\times l)\}\subset\mP^2_x\times(\mP^2_l)^*\times\mP^D_f
\eeq
Therefore
\beq
[\tilde\Sigma]=(L+X)\sum^\frac{p(p+3)}{2}_{i=0}((d-p)X+F)^i(pL)^{\frac{p(p+3)}{2}-i}
\eeq
Extracting the coefficient of $X^2L^2F^{\frac{p(p+3)}{2}-3}$
we get: \beq [\Sigma]=\frac{p(p-1)(p+4)}{8}(d-p)\left(
d(p^2+3p-2)-p(p^2+3p-6)\right)f^{\frac{p(p+3)}{2}-3}\in
H^{p(p+3)-6}(\mP^D_f) \eeq which gives:
\begin{itemize}
\item{$p=2$}: $[\Sigma_{A_2}]=12(d-1)(d-2)$
\item{$p=3$}: $[\Sigma_{E_6}]=21(d-3)(4d-9)$
\item{$p=4$}: $[\Sigma_{W_{12}}]=24(d-4)(13d-44)$
\end{itemize}

\paragraph{The series $x^p_2+x^{p+2}_1~~~(A_3,~E_8,~W_{1,0})$}
The analysis of Newton diagram, similar to the $A_k$ case, shows that the defining conditions are:
\beq\ber
f^{(p+1)}\wedge\underbrace{l\wedge\dots\wedge l}_{p+1}=0\\
f^{(k)}\wedge\underbrace{l\wedge\dots\wedge l}_{k-m}=0,~~\mbox{for}~ k<p+2-\frac{2m}{p},~~m=1,2\dots,k\\
\eer\eeq
In particular the necessary and sufficient set of conditions is:
\beq\ber
f^{(p+1)}\wedge\underbrace{l\wedge\dots\wedge l}_{1+p-[\frac{p-1}{2}]},~~~~f^{(p)}\wedge l=0,~~~~f^{(p-1)}=0
\eer\eeq
As previously, they are solved in terms of auxiliary symmetric form:
\beq\ber
f^{(p+1)}\sim\mbox{\bf{SYM}}\left((l)^{(1+[\frac{p-1}{2}])}),B^{(p-[\frac{p-1}{2}])}\right)\\
B^{(p-[\frac{p-1}{2}])}(x)\sim\underbrace{(l\times\dots\times l)}_{p-1-[\frac{p-1}{2}]}
\eer\eeq
This gives:
\beq\ber
[\tilde\Sigma]=(L+X)\sum^{N_1}_{i=0}((d-p-1)X+F)^i(([\frac{p-1}{2}]+1)l+B)^{N_1-i}
\sum^{N_2}_{j=0}(B+X)^j((p-1-[\frac{p-1}{2}])l)^{N_2-j}\\
N_1={p+3\choose{2}}-1~~~
N_2={p+1-[\frac{p-1}{2}]\choose{2}}-1
\eer\eeq
In particular:
\begin{itemize}
\item{$p=2$} $[\Sigma_{A_3}]=2(25d^2-96d+84)$
\item{$p=3$}  $[\Sigma_{E_8}]=9(15d-56)(3d-8)$
\item{$p=4$}~~$[\Sigma_{W_{1,0}}]=12(105d^2-868d+1773)$
\end{itemize}

\paragraph{The series $x^p_2+x^{p+3}_1\!\!,~p\geq3~~~(J_{10},~W_{18})\!\!$}The defining conditions are:
\mbox{$f^{(p+2)}\wedge\underbrace{l\wedge\dots\wedge l}_{p+2-[\frac{p-1}{3}]}=0$},  ~~
$f^{(p+1)}\wedge\underbrace{l\wedge\dots\wedge l}_{p+1-[\frac{2(p-1)}{3}]}=0,~~f^{(p)}\wedge l=0$
In terms of auxiliary symmetric tensors we have:
\beq\ber
f^{(p+2)}\sim\mbox{\bf{SYM}}\left((l)^{(1+[\frac{p-1}{3}])}),B_1^{(p+1-[\frac{p-1}{3}])}\right)\\
B_1^{(p+1-[\frac{p-1}{3}])}(x)\sim{\bf{SYM}}\left((l)^{([\frac{2(p-1)}{3}]-[\frac{p-1}{3}])},B_2^{(p-[\frac{2(p-1)}{3}])}\right)\\
B_2^{(p-[\frac{2(p-1)}{3}])}(x)\sim(l)^{(p-1-[\frac{2(p-1)}{3}])}
\eer\eeq
This gives: \beq\ber
[\tilde\Sigma]=(L+X)\sum^{N_1}_{i=0}((d-p-2)X+F)^i(((1+[\frac{p-1}{3}]))L+B_1)^{N_1-i}
\sum^{N_2}_{j=0}(B_1+X)^j(([\frac{2(p-1)}{3}]-[\frac{p-1}{3}])L+B_2)^{N_2-j}\\
\sum^{N_3}_{k=0}(B_2+X)^k((p-1-[\frac{2(p-1)}{3}])L)^{N_3-k}\\
N_1={p+4\choose{2}}-1~~~
N_2={p+2-[\frac{p-1}{3}]\choose{2}}-1~~~
N_3={p+1-[\frac{2(p-1)}{3}]\choose{2}}-1
\eer\eeq
In particular:
\begin{itemize}
\item{$p=3$}:  $[\Sigma_{J_{10}}]=2(385d^2-2715d+4611)$
\item{$p=4$}:  $[\Sigma_{W_{18}}]=3(1173d^2-11288d+26608)$
\end{itemize}

\paragraph{The series $x^p_2+x^{p+4}_1,~p\geq4,~(X_{2,0})$}
$f^{(p+3)}\wedge\underbrace{l\wedge\dots\wedge l}_{p+3-[\frac{p-1}{4}]}=0$,
~~$f^{(p+2)}\wedge\underbrace{l\wedge\dots\wedge l}_{p+2-[\frac{(p-1)}{2}]}=0$,~~
$f^{(p+1)}\wedge\underbrace{l\wedge\dots\wedge l}_{p+1-[\frac{3(p-1)}{4}]}=0,~~~f^{(p)}\wedge l=0$
\paragraph{The series $x^p_2x_1+x^q_1,~\frac{p}{2}<q-1\leq2p$}
From the Newton diagram one gets conditions:
\beq\label{DKfirst}
\di^i_2\di^j_1f=0,~\mbox{for}~j<q+\frac{i(1-q)}{p}
\eeq
In several lowest cases we have:
\begin{itemize}
\item{$x^p_2x_1+x^{p-1}_1$}. $J_{10}$ (has been treated previously),~~ $[\Sigma_{W_{17}}]=2584d^2-23664d+53259$
\item{$x^p_2x_1+x^p_1,~A_3,~E_7,~W_{13}$}
Here the conditions (\ref{DKfirst}) mean: $f^{(p)}\wedge l=0,~~f^{(p+1)}\wedge\underbrace{l\wedge\dots\wedge l}_{p+1}=0$, that is:
$f^{(p+1)}=\mbox{\bf{SYM}}(l,B^{(p)}),~~B^{(p)}(x)=(l)^{(p-1)}$. So the cohomology class is:
\beq\ber
(L+X)\sum^{N_1}_{i=0}((d-p-1)X+F)^i(L+B)^{N_1-i}\sum^{N_2}_{j=0}(B+X)^j((p-1)l)^{N_2-j}\\
N_1={p+3\choose{2}}-1~~~N_2={p+1\choose{2}}-1
\eer\eeq
In particular
\begin{itemize}
\item{$p=3$}:  $[\Sigma_{E_7}]=3(84d^2-488d+693)$
\item{$p=4$}:  $[\Sigma_{W_{13}}]=3(273d^2-2132d+4128)$
\end{itemize}

\item{$x^p_2x_1+x^{p+1}_1,~~~D_4,~X_9$}
This defines singularity of the type: ordinary multiple point, and was considered previously.
\item{$x^p_2x_1+x^{p+2}_1,p\geq2~~D_5,~Z_{11}$}.
Here the conditions (\ref{DKfirst}) mean: $f^{(p)}=0,~~f^{(p+1)}\wedge l\wedge l=0$, that is:
$f^{(p+1)}=\mbox{\bf{SYM}}((l)^{(p)},B^{(1)}),~~~B(x)=0$ So the
cohomology class is:
\beq
(L+X)(B+X)\sum^N_{i=0}((d-p-1)X+F)^i(pL+B)^{N-i}~~~N={p+3\choose{2}}-1~
\eeq
In particular
\begin{itemize}
\item{$p=2$}:   $[\Sigma_{D_5}]=12(d-2)(7d-19)$
\item{$p=3$}:  $[\Sigma_{Z_{11}}]=18(267 - 154d + 22d^2)$
\end{itemize}
\item{$x^p_2x_1+x^{p+3}_1,p\geq2~~D_6,~~Z_{13}$}

Here the conditions (\ref{DKfirst}) mean: $f^{(p)}=0,~~f^{(p+1)}\wedge l\wedge l=0,~~~
f^{(p+2)}\wedge\underbrace{l\wedge\dots\wedge l}_{p+3-[\frac{p}{2}]}=0$, that is:
$f^{(p+2)}=\mbox{\bf{SYM}}((l)^{([\frac{p}{2}])},B_1^{(p+2-[\frac{p}{2}])}),~~~
B_1^{(p+2-[\frac{p}{2}])}(x)=\mbox{\bf{SYM}}((l)^{(p-[\frac{p}{2}])},B_2^{(1)}),~~~B_2^{(1)}(x)=0$

So the cohomology class is:
\beq\ber
(L+X)(B_2+X)\sum^{N_1}_{i=0}((d-p-2)X+F)^i([\frac{p}{2}]L+B_1)^{N_1-i}
\sum^{N_2}_{j=0}(B_1+X)^j((p-[\frac{p}{2}])L+B_2)^{N_2-j}\\
N_1={p+4\choose{2}}-1~~~N_2={p+3-[\frac{p}{2}]\choose{2}}-1
\eer\eeq
\begin{itemize}
\item{$p=2$}:~~$[\Sigma_{D_6}]=14(16d^2-87d+114)$
\item{$p=3$}:~~$[\Sigma_{Z_{13}}]=3(286d^2-2148d+3977)$
\end{itemize}
\end{itemize}

\paragraph{Singularities defined by two lines on Newton diagram}Examples:
$X_{1,p}$,~\mbox{$0<p\leq2$}:   $x^4_2+x^2_2x^2_1+x^{4+p}_1~~~~W_{1,p},~0<p\leq2~~x^4_2+x^2_2x^3_1+x^{6+p}_1$

Up to now we have considered quasihomogeneous singularities (defined by vanishing of derivatives corresponding to the points below certain line on
Newton diagram). For higher singularities however the typical case is semi-quasi-homogeneous (with two or more lines on Newton diagram).
In this case the calculations go along the
similar lines, as an example we the cases \mbox{$X_{1,p},~p\leq2$},~~~~\mbox{$W_{1,p},~p\leq2$}. To avoid cumbersome formulae we just
record the defining conditions.
\begin{itemize}
\item{$X_{1,1}$} $f^{(4)}\sim\mbox{\bf{SYM}}(l,l,B^{(2)})~~B^{(2)}(x)=0$
\item{$X_{1,2}$} $f^{(5)}\sim\mbox{\bf{SYM}}(l,B^{(4)})~~B^{(4)}(x)\sim\mbox{\bf{SYM}}(l,C^{(2)})~~C^{(2)}(x)=0$
\item{$W_{1,1}$} $f^{(6)}\sim\mbox{\bf{SYM}}(l,B^{(5)})~~B^{(5)}(x)\sim\mbox{\bf{SYM}}(l,C^{(3)})~~C^{(3)}(x)\sim l\times l$
\item{$W_{1,2}$} $f^{(7)}\sim\mbox{\bf{SYM}}(l,B^{(6)})~~B^{(6)}(x)\sim\mbox{\bf{SYM}}(l,C^{(4)})~~C^{(4)}(x,x)\sim l\times l$
\end{itemize}
\subsubsection{Non-linear strata}\label{SecExplicitCalculationNonLinearStrata}
 As is seen in the case of linear strata, from the enumerative point of view, different
 types of singularities naturally appear in series. For example, for linear singularities: $f(x_1,x_2)=x^p_2+x^p_1+\dots$
 (multiple points), $f(x_1,x_2)=x^p_2+x^{p+1}_1+\dots$
 ($A_2,D_5,E_6$) etc. Similarly the non-linear strata appear in series, all the members of the series have the same "complexity
 of enumeration". We consider here several simplest cases.

 {\bf The "first" nonlinear series: $f(x_1,x_2)=x^p_2+x^{2p+1}_1,~A_4~,E_{12},~W_{24}$ }
The defining (nonlinear) conditions here, mean that the quasi-homogeneous form:
$a_{0,p}x^p_2+a_{2,p-1}x^{p-1}_2x^2_1+\dots+a_{2p,0}x^{2p}_1$ is maximally
degenerated, i.e. it can be written as: $(\alpha x_2+\beta x^2_1)^p$. Therefore,
demanding $a_{0,p}=0~(\alpha=0)$ leaves $x^{2p}_1$ only, i.e.
we obtain the germ: $f(x_1,x_2)=x^p_2x_1+x^{2p}_1$. This stratum arises with multiplicity: $2^{p-1}$.
So in terms of cohomology classes one has:
\beq\label{FirstNonlinear}
[f(x_1,x_2)=x^p_2+x^{2p+1}_1]\left([a_{0,p}=0]-[l_1^p=0]\right)=2^{p-1}[f(x_1,x_2)=x^p_2x_1+x^{2p}_1]
\eeq
The r.h.s. is a linear stratum. Working as in the previous section, its cohomology class is obtained from:
\beq
[l(x)=0][f^{2p-1}\sim\mbox{SYM}(l,l,B^{(2p-3)}_1)]\Pi_{i=1}^{p-2}[B^{(2(p-i)-1)}_i(x)\sim\mbox{SYM}(l,B^{(2(p-i)-3)}_{i+1})]
[B^{(1)}_{p-1})(x)=0]
\eeq
From here one extracts the coefficient of maximal (nonzero) powers of $B_i$. Substitution to (\ref{FirstNonlinear}) gives the answer.
(The general formula is absolutely cumbersome, the answers for cases of $A_4,E_{12},W_{24}$ are written at the beginning of the paper.)

{\bf The "second" nonlinear series: $f(x_1,x_2)=x^p_2+x_2x^{2p-1}_1,~~A_5,~E_{13},~W_{25}$}
As in the previous case, we have maximal degeneracy of the form
 $a_{0,p}x^p_2+a_{2,p-1}x^{p-1}_2x^2_1+\dots+a_{2p,0}x^{2p}_1$. The additional
condition is that the form: $a_{1,p}x^p_2x_1+a_{3,p-1}x^{p-1}_2x^3_1+\dots+a_{2p,0}x^{2p+1}_1$ has a common
root with the maximally degenerated
form. Therefore the degeneration here is(in terms of classes):
\beq
[f(x_1,x_2)\!=\!x^p_2+x_2x^{2p-1}_1]\left([a_{0,p}\!=\!0]\!-\![l_1^p\!=\!0]\right)\!=
\!2^{p-1}[f(x_1,x_2)\!=\!x^p_2x_1+x^{2p+1}_1]+2^{p-1}[f(x_1,x_2)\!=\!x^{p+1}_2+x^{p-1}_2x^3_1+x^{2p}_1]
\eeq
As in the case of the "first" nonlinear series the singularities on the left are linear, their classes are obtained immediately, and one
gets the class of the "second" nonlinear series. (The answers for $p=2,3,4$ are given at the beginning of the text).

{\bf The "third" nonlinear series\footnote{For $p=2$ we get the same $A_5$ as in the "second" nonlinear series}:
$x^p_2+x^{2p+2}_1,~p>2,~E_{14},~W_{26}$}

Here the degeneration is:
\beq
[f(x_1,x_2)\!=\!x^p_2+x^{2p+2}_1]\left([a_{0,p}\!=\!0]-[l_1^p\!=\!0]\right)\!=\!2^{p-1}[f(x_1,x_2)\!=
\!x^p_2x_1+\dots+x^2_2x^{2p-3}_1+x^{2p+2}_1]+2^{p-1}[f(x_1,x_2)\!=\!x^{p+1}_2+x^{2p}_1]
\eeq
Now, one stratum on the r.h.s. is nonlinear, we need additional degeneration to obtain the answer.
We omit the calculations.
\subsection{Tangency of two curves}
Here we consider two conics with prescribed degrees of tangency.
 Start from the local equations:
\beq\ber
C_1(x,y)=ey+ax^2+bxy+cy^2~~~C_2(x,y)=Ey+Ax^2+Bxy+Cy^2
\eer\eeq
The conditions of triple tangency read: $Ea=eA$, $Ab=aB$ and $Eb=eB$.
Therefore we define the stratum of maximally tangent conics as
\beq
\tSi_{2,2}=\{(x,l,C_1,C_2)|l(x)=0,~~rk\bpm C_1\\C_2\\l\times l\epm<3,~C_1^{(1)}|_x\sim l\sim C_2^{(1)}\}\subset
\mP^2_x\times(\mP^2_l)^*\times\mP^5_{C_1}\times\mP^5_{C_2}
\eeq
Note that here we cannot choose mutually transversal sets of conditions. E.g. choosing only $C_1^{(1)}|_x\sim l$
produces a residual piece over $C_1^{(2)}\sim l\times l$ (which is of bigger dimension than the stratum). Similarly,
choosing $C_1^{(1)}|_x\sim g^{(1)}$ produces the residual piece over $C_1\sim C_2$.

We degenerate by $a\ra0$ (the cohomology class of the degeneration is $C_1-2X+2L$), it means either
\li $A=0,~Eb=eB$ (two nodal conics with a common tangent)
\beq
\tSi=\{(l,l_1,l_2,C_1,C_2)|rk\bpm l_1\\l_2\\l\epm<3,~C_1\sim l\times l_1,~~~C_2\sim l\times l_2\}
\eeq
or
\li
$b=0=e$ (a conic tangent to a double line)
\beq
\tSi=\{(x,l,C_1,C_2)|l(x)=0,~~C_1^{(2)}\sim l\times l,~~~C_2^{(1)}\sim l\}
\eeq
 The calculation gives:
\beq\label{ConicMaxTangent}
\ber
[\tilde\Sigma_{2,2}(x,l,C_1,C_2)]=(l+x)\Big(\ber
C_1^4C_2^2+C_1^3C_2^3+C_1^2C_2^4+(l+2x)(C_1^4C_2+C_1C_2^4)+\\
(3l+2x)(C_1^3C_2^2+C_1^2C_2^3)+2lx(C_1^4+C_2^4)+2(4l^2+3x^2)(C_1^3C_2+C_1C_2^3)\\
+6(2l^2+x^2)C_1^2C_2^2+4lx^2(C_1^3+C_2^3)+12lx^2(C_1^2C_2+C_1C_2^2)
\eer\Big)
\eer\eeq

{\it Address}: School of Mathematical Sciences, \mbox{Tel Aviv}
University, \mbox{Ramat Aviv}, 69978 \mbox{Tel Aviv}, Israel.

{\it E-mail}: kernerdm@post.tau.ac.il

\begin{thebibliography}{99}\small\addtolength{\parskip}{-8pt}
\bibitem[AR]{AR}  M.Alberich-Carraminana, J.Roe "Enriques diagrams and adjacency of planar curve singularities",
 Canad. J. Math.  57  (2005),  no. 1, 3--16. Preprint arXiv:math.AG/0108223
\bibitem[AGV]{AGV}
V. I.Arnol'd, S. M.Guse\u\i n-Zade, A. N.Varchenko, "Singularities of differentiable maps. Vol. I.
The classification of critical points, caustics and wave fronts."
Monographs in Mathematics, 82. {\it Birkh\"{a}user Boston, Inc., Boston, MA}, 1985.
\bibitem[Al]{Alufi}P.Aluffi "Characteristic classes of discriminants and enumerative geometry", Comm. Algebra
26 (1998), no. 10, 3165--3193.
\bibitem[CH1]{CapHar1} L. Caporaso, J. Harris: "Enumerating rational curves: the rational fibration method". {\it Compositio
 Math.} {\bf 113} (1998), no. 2, 209--236.
\bibitem[CH2]{CapHar2} L. Caporaso, J. Harris: "Counting plane curves of any genus". {\it Invent. Math.}
{\bf 131} (1998), no. 2, 345--392.
\bibitem[Fr.Itz]{diFrItz} P. di Francesco, C. Itzykson: "Quantum intersection rings", in:
{\it The Moduli Space of Curves} (R. Dijkgraaf, C. Faber, and G. van der Geer, eds.),
Progress in Mathematics, vol. 129, Birkh$\ddot{\rm{a}}$user, Boston, 1995, pp. 81-148
\bibitem[Ful]{Ful} W.Fulton: "Intersection theory." Second edition.
 A Series of Modern Surveys in Mathematics, Springer-Verlag, Berlin, 1998.
\bibitem[G]{Greuel} G.-M.Greuel: "Constant Milnor number implies constant multiplicity for quasihomogeneous singularities".
{\it Manuscripta Math.} {\bf 56 (1986)}, no. 2, 159--166.
\bibitem[GLS1]{Shustin} G.-M. Greuel, C. Lossen, E. Shustin:
"Introduction to Singularities and Deformations", Series: Springer Monographs in Mathematics, 2006, Approx. 400 p.,
\bibitem[GLS2]{GLS} G.-M. Greuel, C. Lossen, E. Shustin: "The variety of plane curves with ordinary singularities
is not irreducible". {\it Internat. Math. Res. Notices} {\bf 11}
(2001), 543--550.
\bibitem[Kaz1]{Kaz1} M. Kazarian. "Classifying spaces of singularities and Thom polynomials". In: {\it New developments in
Singularity Theory (Cambridge 2000)}, NATO Sci. Ser. II Math.
Phys. Chem, 21, Kluwer Acad. Publ.,  Dordrecht, 2001, pp.
117--134.
\bibitem[Kaz2]{Kaz2} M. Kazarian. "Thom polynomials for Lagrange, Legendre, and critical point function
singularities". {\it Proc. Lond. Math. Soc.} (3) {\bf 86} (2003),
707--734.
\bibitem[Kaz3]{Kaz3} M. Kazarian. "Multisingularities, cobordisms, and enumerative geometry".
{\it Russ. Math. Surveys} {\bf 58(4)} (2003), 665--724.
\bibitem[Kaz4]{Kaz4} M. Kazarian. "Characteristic Classes in Singularity theory",
Doctoral Dissertation (habilitation thesis), {\it Steklov Math. Inst.},2003.
\bibitem[Ker]{Ker}D.Kerner {\it Enumeration of singular algebraic curves}, Israel
Journal of Math. 155 (2006), pp1-56, arXiv:math.AG/0407358
\bibitem[Kl1]{Klei1} S.L. Kleiman, "The enumerative theory of singularities". Real and complex singularities
({\it Proc. Ninth Nordic Summer School/NAVF Sympos. Math., Oslo, 1976}), pp. 297--396. Sijthoff and Noordhoff, Alphen aan den Rijn, 1977.
\bibitem[Kl2]{Klei2} S.L. Kleiman, "Intersection theory and enumerative geometry: a decade in review".
With the collaboration of Anders Thorup on §3. {\it Proc. Sympos. Pure Math., 46, Part 2,
Algebraic geometry, Bowdoin, 1985 (Brunswick, Maine, 1985), 321--370, Amer. Math. Soc., Providence, RI, 1987.}
\bibitem[KlPi1]{KleiPien1} S.Kleiman, R.Piene, "Enumerating singular curves on surfaces".{\it  Algebraic geometry: Hirzebruch 70
(Warsaw, 1998), 209--238, Contemp. Math., 241, Amer. Math. Soc., Providence, RI, 1999}.)
\bibitem[KlPi2]{KleiPien2} S.Kleiman, R.Piene, "Node polynomials for families: methods and applications".
 {\it Math. Nachr. 271, 69--90 (2004)}.
\bibitem[Li]{Liu} A.-K. Liu. "The Algebraic Proof of the Universality Theorem", Preprint
arXiv:math.AG/0402045.
\bibitem[Lu]{Luengo} I.Luengo,  "The $\mu$-constant stratum is not smooth". {\it Invent. Math.} {\bf 90} (1987), no. 1, 139--152
\bibitem[PH]{Pham}F.Pham, "Remarque sur l'\'{e}quisingularit\'{e} universelle", Pr\'epublication Universit\'e de Nice
Facult\'{e} des Sciences, 1970.
\bibitem[R]{Ran}Z. Ran. "Enumerative geometry of singular plane curves". {\it Invent. Math.} {\bf 97} (1989), no. 3, 447--465.
\bibitem[V1]{Vain1}I. Vainsencher. "Counting divisors with prescribed singularities". {\it Trans. Amer. Math. Soc.}
{\bf 267} (1981), no. 2, 399--422.
\bibitem[V2]{Vain2} I. Vainsencher, "Hypersurfaces with up to six double points".
{\it Special issue in honor of Steven L. Kleiman. Comm. Algebra 31 (2003), no. 8, 4107--4129.}
\end{thebibliography}
\end{document}